\documentclass[a4paper,11pt]{article}
\usepackage{amssymb,amsmath,latexsym,amsthm}

\usepackage{xcolor}
\usepackage{url}
\usepackage{enumerate}
\usepackage[square,sort,comma,numbers]{natbib}
\usepackage{dsfont, esint}
\usepackage{diagbox}
 \usepackage{geometry}
\usepackage{graphicx}
\usepackage{colortbl}
\usepackage{bbm}

\usepackage{endnotes}

\renewcommand{\textbf}[1]{\begingroup\bfseries\mathversion{bold}#1\endgroup}

\usepackage{breqn}

\usepackage{fancyhdr}
\usepackage{pstricks,pst-plot,pst-node,pstricks-add}

\newtheorem{thm}{Theorem}[section]

\newtheorem{corollary}[thm]{Corollary}
\newtheorem{prop}[thm]{Proposition}

\newtheorem{lemma}[thm]{Lemma}

\theoremstyle{definition}
\newtheorem{remark}[thm]{Remark}

\newtheorem{example}[thm]{Example}

\newcommand{\R}{\mathbb R}

\newcommand{\N}{\mathbb N}

\numberwithin{equation}{section}

\def\XXint#1#2#3{{\setbox0=\hbox{$#1{#2#3}{\int}$}
    \vcenter{\hbox{$#2#3$}}\kern-.5\wd0}}

\allowdisplaybreaks

\makeatletter
\def\blfootnote{\xdef\@thefnmark{}\@footnotetext}
\makeatother

\date{date}
\begin{document}

\title{Almost sure recovery in quasi-periodic structures}

\author{Mircea Petrache\footnote{PUC, \texttt{mpetrache@mat.uc.cl}. ORCID id: 0000-0003-2181-169X} and 
Rodolfo Viera\footnote{PUC, \texttt{rodolfo.viera@mat.uc.cl}. ORCID id: 0000-0002-0810-6374}}\blfootnote{Pontificia Universidad Católica de Chile, Facultad de Matematicas, Avda. Vicuña Mackenna 4860, Santiago,}\blfootnote{6904441, Chile}

\date{December 22, 2021}
\maketitle
\begin{abstract}
    We study random perturbations of quasi-periodic uniformly discrete sets in the $d$-dimensional euclidean space. By means of Diffraction Theory, we find conditions under which a quasi-periodic set $X$ can be almost surely recovered from its random perturbations. This extends the recent periodic case result of Yakir \cite{Yakir}.
\end{abstract}

\section*{Introduction}\label{ssection: intro} 
 
Random perturbations of discrete structures have been studied by many authors from different viewpoints. For instance, a problem of interest in statistical physics is to study what happens with the statistical properties (e.g, two-point correlation, structure factor) of the perturbed discrete sets; we refer to \cite{gab,kkt} and the references therein. Another problem of interest in Diffraction theory is to compare the diffraction measure of a random perturbation with the diffraction measure of the deterministic discrete set which is perturbed; see for instance \cite[Chapter 11]{BG1} and \cite{Hof} for more details. Very recently, in \cite{Yakir} Yakir deals with random perturbations of a lattice where it is shown that, under some assumptions on the random displacements, the lattice can be recovered almost surely from the perturbations. In this work we focus on random perturbations and reconstructibility in quasi-periodic structures, including quasicrystals, and on weakening the independence conditions on the random perturbations. We stress that quasicrystals exist not only as mathematical objects, but also as naturally ocurring \cite{quasicrystal2015} physical materials, which although rare, have been argued to be ubiquitous throughout the universe \cite{quasicrystal2020}. Since the discovery of quasicrystalline structures in materials by Shechtman and his collaborators \cite{quasi1} many efforts have been dedicated to understanding their underlying order and their stability under perturbations (see for instance \cite{Hof, Hof2, hypquas} just to mention a few). 
\medskip

In more precise terms, if $X\subset\R^d$ is a quasi-periodic set, we let $(\xi_p)_{p\in X}$ be a sequence of identically distributed random vectors defined in a probability space $(\Omega, \mathcal{F},\mathbb{P})$, with common distribution $\xi$. A {\em random perturbation} of $X$ is a realization of the random set
 
 \begin{equation*}
     X_{\xi}:=\{p+\xi_p(\omega):\ p\in X\}.
 \end{equation*}
 
For two random processes $X,Y\subset\mathbb{R}^d$ that take values into locally finite configurations, we say that $X$ is {\em almost surely recovered} from $Y$ if there exists a measurable function $F$ defined over the essential range of the process $Y$ such that almost surely $F(Y)=X$; this notion extends naturally to random measures. If $X$ is a fixed (deterministic) subset of $\R^d$, we identify $X$ with the point process which takes constant value equal to $X$. \\

In this work we focus on finding conditions under which the Fourier transform of $X_{\xi}$ determines the Fourier transform of a quasi-periodic set $X$ almost surely, based on an explicit formula that links these two Fourier transforms. It follows that $X$ can then be almost surely recovered from $X_\xi$, but we can say even more. In fact, some mathematical and physical point configurations can be classified in terms of their hidden order (e.g, repetitivity of separated nets \cite{LagPle}, quasi-periodicity \cite{BG1, apercrys}, Hyperuniformity \cite{Torstil, hypquas}). In this context, our main results (Theorems A,B and C below) provide another way to distinguish the large class of point sets that have diffraction measure singular with respect to Lebesgue measure (see hypotheses (H1), (H2) in Section \ref{ssec:hypo}) in terms of the robustness of their Fourier Transform. This allows to recognize these configurations under random perturbations, by taking a Fourier Transform. Note that as a special case our setting includes quasicrystals (defined by asking that their diffraction measure is purely atomic), which in turn include the classical case of lattices, for which the diffraction measure is supported on their dual lattice. In \cite{BF} the authors show that SCD patterns have singular diffraction measure (in particular they satisfy our hypotheses), but are not quasicrystals.

\medskip

{\bf Notations. }We include here some basic definitions and notations which we will use throughout this work. We use the notation $|\cdot|$ indistinguishably for the euclidean norm in $\mathbb{R}^d$ and for the module (resp. absolute value) of a complex number (resp. real number), which will be clear in each context. As usual we denote by $B(x,R):=\{y\in\mathbb{R}^d: |x-y|<R\}$ the euclidean open ball centered at $x\in\mathbb{R}^d$ with radius $r>0$ and we use the notation $B_R:=B(0,R)$.

\medskip

We denote by $Vol$ the Lebesgue measure in $\mathbb{R}^d$ and we denote the cardinality of a finite set $A\subset\mathbb{R}^d$ by $\# A$. Observe that a locally integrable function $f\in L_{loc}^1(\mathbb{R}^d)$ can be identified with $f\cdot Vol$; see section \ref{ssec:diffraction} for more details. Furthermore, for an atomic measure $\mu$ we write $\mu(p)$ instead of $\mu(\{p\})$. Finally,  for strictly positive real-valued functions $f,g$ we use the standard Landau asymptotic notations $O(f(R))$ and $o(g(R))$.

\medskip

{\bf Acknowledgements: }\emph{ M. P. was sponsored by the Chilean Fondecyt Regular grant number 1210462 entitled ``Rigidity, stability and uniformity for large point configurations'', and R. V. was sponsored by the Chilean Fondecyt Postdoctoral grant number 3210109 entitled ``Geometric and Analytical aspects of Discrete Structures''.}
\section{Hypotheses on the deterministic set $X$}\label{ssec:hypo}

%We say that a point set $X\subset\mathbb{R}^d$ is {\em uniformly discrete} if there exists $r>0$ such that $d(x,y)\geq r$ for every $x,y\in X$; in addition, $X$ is called {\em relatively dense} if there is a positive number $\rho$ such that for every $z\in\mathbb{R}^d$ there is $x\in X$ for which $d(x,z)\leq\rho$. Throughout this work we assume that $X$ is a {\em Delone set}, which means that $X$ is relatively dense and uniformly discrete.
%
%\medskip

Throughout this work we assume that $X$ is an {\em uniformly discrete set}, which means that there exists $r>0$ such that $||x-y||\geq r$ for every $x,y\in X$. %\rov{We point out we never use the relatively denseness of $X$. Moreover, there are examples of uniformly discrete sets which are not Delone sets (e.g the set of visible points of $\mathbb{Z}^d$). The principle is that the diffraction spectrum of $X$ does not change by adding or deleting points of zero-density.}
The set $X$ is assumed to satisfy the following two conditions:

\begin{enumerate}
    \item[(H.1)] There exists a positive constant $\mathrm{dens}(X)$, called the {\em asymptotic density} of $X$, such that for every $R>0$ we have
    
    \begin{equation*}
        |\# (X\cap B_R)-\mathrm{dens}(X)\cdot Vol(B_R)|=o(R^d).
    \end{equation*}
    
%Equivalently, the point set $X$ has an {\em asymptotic density}, which means that the following limit exists:
    
 %   \begin{equation*}\label{asympdens}
  %      \mbox{dens}(X)=\lim_{R\to\infty}\dfrac{\# (X\cap B_R)}{Vol(B_R)}.
   % \end{equation*}
    
    \item[(H.2)] The  diffraction measure $\widehat{\gamma_X}$ of $\delta_X:=\sum_{p\in X}\delta_p$, given by
    \begin{equation}\label{weaklim}
    \forall f\in C^\infty_c(\mathbb R^d),\qquad \langle \widehat{\gamma_X}, f\rangle:=\lim_{R\to\infty}\frac{1}{Vol(B_R)}\int_{\mathbb R^d} f(\lambda)\left|\sum_{p\in X\cap B_R} e^{-2\pi i\langle p, \lambda\rangle}\right|^2 d\lambda
    \end{equation}
    exists and is singular with respect to Lebesgue measure, where $\delta_p$ denotes the Dirac mass at the point $p\in X$. 
\end{enumerate}

There are plenty of uniformly discrete sets $X$ satisfying the hypothesis (H.1) (where trivially all lattices are included). Non-trivial examples are cut-and-project sets constructed with an appropriate choice of the window; see Section \ref{exstat} for precise definitions. A particular case are the cut-and-project sets which are {\em bounded displacement equivalent} to a lattice; these sets include, for instance, the set of vertices of the {\em Penrose tilling}. It is worth mentioning that a generic cut-and-project set in $\mathbb{R}^d$ is, up to bounded displacement, equivalent to a lattice (for more details see \cite{HKW}). Another interesting example of a set $X$ satisfying (H.1) is the $X=V_d$, known as the set of {\em visible points} of $\mathbb{Z}^d$, and defined as $V_d:=\{x\in\mathbb{Z}^d:\ gcd(x)=1\}$, where $gcd(x)$ denotes the greatest common divisor of the coordinates of $x\in\mathbb{Z}^d$; it is known that $\mathsf{dens}(V_d)=1/\zeta(d)$, where $\zeta$ denotes the {\em Riemann's zeta function}. We refer to \cite{BMP} for the proof of (H.1) in this case. Observe that $V_d$ is not relatively dense, which roughly speaking amounts to saying that $V_d$ has arbitrarily large holes. % \rov{Furthermore, in \cite{RSW} it was proven that the set of irreducible cut-and-project sets satisfying (H.1) is a set of full-measure, with an appropriate measure defined in the space of irreducible cut-and-project sets (namely, the so-called {\em Ratner-Marklof-Str\"{o}mbergsson measures}). Borramos esto \'ultimo?}. 

\medskip

All the abovementioned examples satisfy (H.2). See for instance Corollary 9.3 and Proposition 9.9 in \cite{BG1} and  \cite[Section 9.4]{BG1} for several examples of cut-and-project sets and the precise computations of their diffraction measures; see also \cite{BMP} for a proof that $\widehat{\gamma_{V_d}}$ is purely atomic in the case of visible points of a lattice. In particular a set $X$ satisfies (H.2) if its weak dual $X^*$ is discrete, where $X^*$ is the support of the Fourier transform of $\delta_X:=\sum_{p\in X}\delta_p$, given by $\widehat{\delta_X}=\sum_{\lambda\in X^*}A(\lambda)\delta_\lambda$. More explicitly, there exist $\{A(\lambda)\}_{\lambda\in X^*}$ such that
    \begin{equation*}
        \displaystyle\lim_{R\to\infty}\frac{1}{Vol(B_R)}\sum_{p\in X\cap B_R}e^{-2\pi i\langle p,\lambda\rangle}=\left\{\begin{array}{ll} A(\lambda) &\text{ for }\lambda\in X^{*},\\
   0&\text{ else.}
   \end{array}\right.
   \end{equation*}
 Note that $X^*$ coincides with the usual dual for lattices and cut-and-project sets. %See Section \ref{exstat} for a example of a concrete Delone set $X$ satisfying (H.1)-(H.3).
 
 \medskip
 
 We consider the general setting given by (H.2) since there are uniformly discrete point configurations in $\mathbb{R}^d$ for which only the absolutely continuous part of their diffraction spectrum vanishes; see for instance \cite{BF}.
%%%%%%%%%%%%%%%%%%%%%%%%%%%%%%%%%%%%%%%%%%%%%%%%%%%%%%%%%%%%%%%%%%%%%%%%%%%%%%%%%%%%%%%%%%%%%%%%%%%%%%%%%%%%%%%%%%%%%%%%%%%%%%%%%%%%%%%%%%%%%%%%%%%%%%%%%%%%%%%%%%%%%%%%%%%%%%%%%%%%%%%%%%%%%%%%%%%%%%
\section{Statements of the results}

 Recall that $X_{\xi}$ is the realization random set $\{p+\xi_p:\ p\in X\}$. For every $\lambda\in\mathbb{R}^d$, $\omega\in\Omega$ and $R>0$, we define $M_{R,\xi,\omega}(\lambda)$ by

\begin{equation*}
    M_{R,\xi,\omega}(\lambda):=\dfrac{1}{Vol(B_R)}\sum_{p\in X_{\xi}\cap B_R}e^{-2\pi i\langle p,\lambda\rangle}.
\end{equation*}

The aim of this work is to show that under our hypotheses (H.1) and (H.2) the Fourier transform %diffraction spectrum
of $X$ is preserved under random perturbations. As a byproduct, this allows us to recover almost surely the quasi-periodic set $X$ from $X_{\xi}$.\\

{\bf Theorem A. }\label{Main Thm}{\em Suppose $X\subset \R^d$ satisfies hypotheses (H.1) and (H.2) from Section \ref{ssec:hypo}. Assume that the $\xi_p$'s are i.i.d and there exists a positive number $\varepsilon$ such that $\mathbb{E}(|\xi|^{d+\varepsilon})<\infty$. Then, almost surely the following holds for every $\lambda\in\mathbb{R}^d$,in the sense of vague convergence of measures:}

\begin{equation}\label{eqthmA}
    \displaystyle\lim_{R\to\infty}M_{R,\xi,\omega}(\lambda)=\mathbb{E}\left[e^{-2\pi i\langle\xi,\lambda\rangle}\right]\widehat{\delta_X}(\lambda),
\end{equation}

{\em where $\widehat{\delta_X}$ is the (weak-)Fourier transform of $X$ given by $\lim_{R\to\infty}\frac{\widehat{\delta_{X\cap B_R}}}{Vol(B_R)}$.}
\medskip

A key step in the proof of Theorem A is the strong law of large numbers. Since the hypothesis of this well-known result can be weakened, we can adapt our proof to show that the conclusions of Theorem A still hold for non-independent identically distributed random vectors under mild assumptions on the correlations (see section \ref{ssec: non-ind rv} for the precise definitions).\\

\noindent{\bf Theorem B. }\label{Thm B}{\em Assume that $X$ satisfies (H.1) and (H.2), that the $\xi(p), p\in X$ are identically distributed, and that there is $\varepsilon>0$ such that $\mathbb{E}(|\xi|^{d+\varepsilon})<\infty$. If $(\xi_p)_{p\in X}$ is a compactly-mixing random process such that for $Y_p:=|\xi_p|\mathbbm{1}_{|\xi_p|\le |p|}$ there holds
\begin{equation}\label{covhyp}
    \sum_{N\ge 1}\frac{1}{N^2}\sum_{\substack{p,q\in B_N\\ p\neq q}}|\mathsf{Cov}(Y_p,Y_q)|<\infty,
\end{equation}
then almost surely the limit (\ref{eqthmA}) still holds for every $\lambda\in\mathbb{R}^d$.
}
\medskip

We note that a stronger version of Theorems A, B is achievable if we assume that $X$ is an asymptotically affine perturbation of a lattice. In this case we only need the $\xi_p$ to be weakly mixing. See Theorem D and Section \ref{asyaffine} for details.

\medskip
Though the conclusion of Theorem B holds for every uniformly discrete set $X$ satisfying hypotheses (H.1) and (H.2), we require strong assumptions on the correlations. However, if the random process $(\xi_p)_{p\in X}$ is such that the random perturbation $X_{\xi}$ is stationary, we have the same conclusion of Theorem B but under more classical conditions on the correlations, namely weak-mixing.

\medskip

\noindent{\bf Theorem C. }{\em Assume that $X$ satisfy the hypotheses (H.1) and (H.2) and that $(\xi_p)_{p\in X}$ is a process of random vectors for which there is a positive number $\varepsilon$ such that $\mathbb{E}[|\xi|^{d+\varepsilon}]<\infty$. If the random perturbation $X_{\xi}$ is stationary under the action of $\R^d$ on itself by translations and weakly-mixing, then almost surely, the limit \ref{eqthmA} holds for every $\lambda\in\mathbb{R}^d$.

%\begin{equation*}
 %   \lim_{R\to\infty} M_{R,\xi,\omega}(\lambda)=\lim_{R\to\infty}\frac{1}{Vol(B_R)}\sum_{p\in X\cap B_R}e^{-2\pi i\langle p,\lambda\rangle}\mathbb{E}\left[e^{-2\pi i\langle\xi_p,\lambda\rangle}\right].
%\end{equation*}

%In particular, if the $\xi_p$'s are identically distributed, we get that almost surely, for every $\lambda\in \R^d$
%\begin{equation*}
 %   M_{R,\xi,\omega}(\lambda)\longrightarrow\mathbb{E}\left[e^{-2\pi i\langle\xi,\lambda\rangle}\right]\widehat{\delta_X}(\lambda).
%\end{equation*}
}

\medskip

As we mentioned before, if the limit \eqref{eqthmA} holds and  $\mathbb{E}(e^{-2\pi i\langle\xi,\lambda\rangle})\neq 0$, this allow us to recover the Fourier transform $\widehat{\delta_X}$ from $\widehat{\delta_{X_{\xi}}}$ without requiring specific conditions on $X$. The above condition on $\xi$ holds in several cases of interest, for example for $\xi$ a mixture of Gaussians.

\begin{corollary}
Suppose that either the hypotheses of Theorem A, B or C hold. If for all $\lambda\in\R^d$ there holds $\mathbb{E}(e^{-2\pi i\langle\xi,\lambda\rangle})\neq 0$, then $\widehat{\delta_{X}}$ is almost surely recoverable from $\widehat{\delta_{X_{\xi}}}$. More precisely, almost surely, for all $\lambda\in\mathbb{R}^d$ there holds

\begin{equation*}
    \widehat{\delta_X}(\lambda)=\mathbb{E}\left(e^{-2\pi i\langle\xi,\lambda\rangle}\right)^{-1}\widehat{\delta_{X_{\xi}}}(\lambda).
\end{equation*}
\end{corollary}

A class of results closely related to our work concern the recovery of the structure factor $S_X$ from the one random perturbations $S_{X_\xi}$. In fact, \cite{kkt,gab} obtain a formula for $S_{X_{\xi}}$ in terms of the structure factor $S_X$ of the unperturbed data $X$ (here $X$ itself is assumed to be a stationary random process). This yields a formula related to \eqref{eqthmA} but in expectation form. This relation is discussed in Section \ref{strfact}.

%%%%%%%%%%%%%%%%%%%%%%%%%%%%%%%%%%%%%%%%%%%%%%%%%%%%%%%%%%%%%%%%%%%%%%%%%%%%%%%%%%%%%%%%%%%%%%%%%%%%%%%%%%%%%%%%%%%%%%%%%%%%%%%%%%%%%%%%%%%%%%%%%%%%%%%%%%%%%%%%%%%%%%%%%%%%%%%%%%%%%%%%%%%%%%%%%%%%%%%
\section{Basics on Mathematical Diffraction Theory}\label{ssec:diffraction}
In this section we introduce some terminology and basic concepts which we use along this work. Let $\mathcal{S}(\mathbb{R}^d)$ be the {\em Schwartz space} of rapidly decreasing $C^{\infty}$ (complex-valued) test functions. The {\em Fourier transform} of $f\in\mathcal{S}(\mathbb{R}^d)$ is defined by

\begin{equation*}
    \widehat{f}(y):=\int f(x)e^{-2\pi i\langle x,y\rangle}dx.
\end{equation*}

By a (complex) {\em measure} we mean a linear functional on the set of compactly supported functions $C_c(\mathbb{R}^d)$ such that for every compact set $K\subset\mathbb{R}^d$ there is a positive constant $M_K$ satisfying that

\begin{equation*}
    |\mu(f)|\leq M_K||f||_{\infty},
\end{equation*}

for every $f\in C_c(\mathbb{R}^d)$ supported on $K$, and where $||\cdot||_{\infty}$ denotes the supremum norm on $C_c(\mathbb{R}^d)$. By the Riesz Representation Theorem, there is an equivalence between this definition of measure and the classical measure-theoretic concept of regular Radon measure. When a measure $\mu$ determines a tempered distribution $T_{\mu}(f)=\mu(f)$ (i.e a continuous linear functional defined on $\mathcal{S}(\mathbb{R}^d)$) via the formula

\begin{equation*}
    \forall f\in\mathcal{S}(\mathbb{R}^d),\qquad T_{\mu}(f):=\int f(x)d\mu(x),
\end{equation*}

then its Fourier transform is defined by $\widehat{\mu}(f):=\mu(\widehat{f})$.\\

The {\em total variation} $|\mu|$ of a measure $\mu$ is the smallest positive measure for which $|\mu(f)|\leq |\mu|(f)$ holds for every non-negative $f\in C_c(\mathbb{R}^d)$. A measure $\mu$ is said to be {\em translation bounded} if for every compact $K\subset\mathbb{R}^d$ there is a positive number $a_K$ such that $|\mu|(x+K)\leq a_K$, for every $x\in\mathbb{R}^d$. \\

Let $P_{\mu}:=\{x\in\mathbb{R}^d:\ \mu(\{x\})\neq 0\}$ be the set of atoms of $\mu$. We say that $\mu$ is {\em purely atomic} (or {\em pure point}) if it has atoms only, i.e, $\mu(A)=\sum_{x\in A\cap P_{\mu}}\mu(\{x\})$ for every Borel set $A$; moreover $\mu$ is called {\em continuous} if $\mu(\{x\})=0$ for every $x\in\mathbb{R}^d$. A measure $\mu$ is said to be {\em absolutely continuous} with respect to a measure $\nu$, if there exists $g\in L^{1}_{loc}(\nu)$ such that $\mu=g\nu$. Furthermore $\mu$ and $\nu$ are called {\em mutually singular} (which is denoted by $\mu\perp\nu$) if there exists a relatively compact Borel set $S\subset\mathbb{R}^d$ such that $|\mu|(S)=0$ and $|\nu|(\mathbb{R}^d\setminus S)=0$ (with obvious interpretation for the value of $|\mu|(A)$).\\

 For $f\in\mathcal{S}(\mathbb{R}^d)$, define $\widetilde{f}(x):=\overline{f(-x)}$, and for a measure $\mu$ we consider $\widetilde{\mu}(f):=\overline{\mu(\widetilde{f})}$. We say that $\mu$ is {\em positive definite} if $\mu(f*\widetilde{f})\geq 0$ for every $f\in C_c(\mathbb{R}^d)$. By Proposition 8.6 in \cite{BG1}, every positive definite measure is Fourier transformable, and its Fourier transform is a positive translation bounded measure.\\
 
 Recall that the {\em convolution} $\mu*\nu$ between two measures $\mu$ and $\nu$ is defined by 

\begin{equation*}
    \mu*\nu(f):=\displaystyle\int f(x+y)d\mu(x)d\nu(y);
\end{equation*}

If $\mu$ is translation bounded and $\nu$ is finite, then the convolution $\mu*\nu$ is well-defined and is also a translation bounded measure. Moreover, whenever $\widehat{\mu}$ is a measure, one has the identity $\widehat{\mu*\nu}=\widehat{\mu}\widehat{\nu}$ (see Theorem 8.5 in \cite{BG1}).\\

We say that a net of measures $(\mu_R)_{R>0}$ (resp. positive definite measures) {\em converges vaguely} to a measure $\mu$ if for every test function $f\in C_c(\mathbb{R}^d)$ (resp. for $f\in\mathcal{S}(\mathbb{R}^d)$)

\begin{equation*}
    \lim_{R\to\infty}\langle\mu_R,f\rangle=\langle\mu,f\rangle,
\end{equation*}

where we use the notation $\langle\mu,f \rangle:=\mu(f)$.\\

We denote by $\mu_R$ the restriction of $\mu$ over the ball $B_R$. The {\em autocorrelation measure of $\mu$} is defined by

\begin{equation*}
    \gamma_{\mu}:=\lim_{R\to\infty}\dfrac{\mu_R*\widetilde{\mu_R}}{Vol(B_R)},
\end{equation*}

whenever this limit exists. In this case, $\gamma_{\mu}$ is the limit of positive definite measures, thus is positive definite, and then it is Fourier transformable and its Fourier transform is positive definite and translation bounded. We also call $\widehat{\gamma_{\mu}}$ the {\em diffraction measure of $\mu$}.\\ %(which will be the analogous of the``spectral measure" in \cite{Yakir}).\\

If $\mu:=\sum_{p\in X}\mu(p)\delta_p$, where $\delta_p$ is the Dirac measure at $p\in X$, to deal with the diffraction measure of $\mu$ we need to prove that, in the sense of vague convergence of measures,

\begin{equation*}
    \lim_{R\to \infty} \frac{1}{Vol(B_R)}\left|\sum_{p\in X\cap B_R}\mu(p)e^{-2\pi i \langle p,y\rangle}\right|^2 dy= \widehat{\gamma_{\mu}}
\end{equation*}

In more abstract terms, we are required to prove the following:
\begin{equation}\label{F-lim=lim-F}
    \lim_{R\to\infty}\mathcal F\left(\frac{1}{Vol(B_R)} \mu_R\ast \widetilde{\mu_R}\right) =\mathcal F\left(\lim_{R\to\infty} \frac{1}{Vol(B_R)}\mu_R\ast \widetilde{\mu_R}\right),
\end{equation}

where $\mathcal{F}$ denotes the Fourier transform. This follows from Lemma 4.11.10 in \cite{BG2}, which states that a net of positive definite measures converging weakly in the sense of measures then their Fourier transforms converge weakly as well.\\

%Translation bounded measures provide another way to make more tractable computations of diffraction measures. More precisely, if $\mu$ is translation bounded, then by Lemma 1.2 in \cite{Sch} we have that

%\begin{equation}\label{auto-trbd}
 %  \gamma_{\mu}=\lim_{R\to\infty}\dfrac{\mu_R*\widetilde{\mu}}{Vol(B_R)},
%\end{equation}

%and thus, by (\ref{F-lim=lim-F}) it follows that 

%\begin{equation}\label{calc_diff}
 %  \widehat{\gamma_{\mu}}=\left(\lim_{R\to\infty}\dfrac{\widehat{\mu_R}}{Vol(B_R)}\right)\overline{\widehat{\mu}}.
%\end{equation}

%%%%%%%%%%%%%%%%%%%%%%%%%%%%%%%%%%%%%%%%%%%%%%%%%%%%%%%%%%%%%%%%%%%%%%%%%%%%%%%%%%%%%%%%%%%%%%%%%%%%%%%%%%%%%%%%%%%%%%%%%%%%%%%%%%%%%%%%%%%%%%%%%%%%%%%%%%%%%%%%%%%%%%%%%%%%%%%%%%%%%%%%%%%%%%%%%%%%%%%%%
\section{Proofs}\label{proofs}
The proofs of Theorems A,B and C rely on Lemma \ref{controlled mass} and Proposition \ref{reduction Main Thm} below. Lemma \ref{controlled mass} is a slightly variation of Lemma 2 in \cite{Yakir}, with the very same proof. It claims that the cardinality of points in $X\cap B_R$ which escape from $B_R$ under random perturbations converges in mean to 0 almost surely (the same happens for the cardinality of points in $X\setminus B_R$ with perturbations lying inside to $B_R$).

\begin{lemma}\label{controlled mass}
Assume that $X$ satisfies hypothesis (H.1) and that there is $\varepsilon>0$ such that $\sup_{p\in X}\mathbb{E}(|\xi_p|^{d+\varepsilon})<\infty$. Then almost surely 
 
 \begin{equation*}
 \begin{split}
     \lim_{R\to\infty}\dfrac{\#\{p\in X:\ |p|\leq R, |p+\xi_p(\omega)|>R\}}{Vol(B_R)}=0,
\end{split}
\end{equation*}
and
 \begin{equation*}
     \lim_{R\to\infty}\dfrac{\#\{p\in X:\ |p|>R, |p+\xi_p(\omega)|\leq R\}}{Vol(B_R)}=0
\end{equation*}
\end{lemma}

Proposition \ref{reduction Main Thm} below roughly speaking says that, under the hypotheses of Theorems A, B or C, the limit defining (\ref{eqthmA}) holds when the indices in the sum from $M_{R,\xi,\omega}$ are taken over $X\cap B_R$ instead of $X_{\xi}\cap B_R$. %Along this work our goal is to show Proposition \ref{reduction Main Thm} under the different hypotheses of Theorems A, B and C. 

\begin{prop}\label{reduction Main Thm}
Suppose that either the hypotheses of Theorem A, B or C hold. Then  almost surely for every $\lambda\in\mathbb{R}^d$ we have that 
\begin{equation}\label{yakirlem3}
    \lim_{R\to\infty}\frac{1}{Vol(B_R)}\sum_{p\in X\cap B_R}e^{-2\pi i\langle p, \lambda\rangle}\left(e^{-2\pi i\langle \xi_p,\lambda\rangle} - \mathbb{E}\left[e^{-2\pi i \langle \xi,\lambda\rangle}\right]\right)=0.
\end{equation}
\end{prop}

As in \cite{Yakir}, the conclusions of Theorems A, B and C follow from Lemma \ref{controlled mass} and Proposition \ref{reduction Main Thm} since

\begin{multline*}
    \sup_{\lambda\in\mathbb{R}^d}\left|M_{R,\xi,\omega}(\lambda)-\frac{1}{Vol(B_R)}\sum_{p\in X\cap B_R}e^{-2\pi i\langle p+\xi_p,\lambda\rangle}\right|\\
    \leq\dfrac{\#\{p\in X:|p|\leq R,\ |p+\xi_p|>R\}}{Vol(B_R)}+\dfrac{\#\{p\in X:\ |p|> R,\ |p+\xi_p|\leq R\}}{Vol(B_R)}\longrightarrow 0,
\end{multline*}

where the last convergence holds almost surely by Lemma \ref{controlled mass}. Hence by the above limit and by Proposition \ref{reduction Main Thm} there holds \eqref{eqthmA}.

\medskip

We now direct our attention to proving Proposition \ref{reduction Main Thm} for $X$ satisfying the hypotheses (H.1) and (H.2). To do this in the context of mathematical diffraction theory, we define

\begin{equation}\label{muxlamb}
\mu_X^{\lambda}:=\sum_{p\in X}e^{-2\pi i \langle p, \lambda\rangle}\delta_p
\end{equation} 

and 

\begin{equation}\label{muxilam}
\mu_\xi^{\lambda}:=\sum_{p\in X}\left(e^{-2\pi i\langle \xi_p,\lambda\rangle} - \mathbb{E}[e^{-2\pi i \langle \xi,\lambda\rangle}]\right)\delta_p.
\end{equation}

By extending the strategy of \cite{Yakir}, we will aim to show that the \emph{diffraction measures} of $\mu_X^{\lambda}$ and $\mu_{\xi}^{\lambda}$, denoted by $\widehat{\gamma_\xi^\lambda}:=\widehat{\gamma_{\mu_\xi^\lambda}}$ and $\widehat{\gamma_X^\lambda}:=\widehat{\gamma_{\mu_X^\lambda}}$ respectively, exist and are almost surely mutually singular. In fact, we will prove in Lemma \ref{lem1} that $\widehat{\gamma_{X}^{\lambda}}$ is singular with respect to Lebesgue measure %purely atomic
and in Propositions \ref{prop2}, \ref{noind-prop2} and \ref{stationarycase} below that $\widehat{\gamma_{\xi}^{\lambda}}$ is almost surely absolutely continuous with respect to the Lebesgue measure, under the hypotheses of Theorems A,B and C respectively. Thus the fact that $\widehat{\gamma_X^\lambda}\perp\widehat{\gamma_\xi^\lambda}$ is then used in combination with Theorem \ref{yakirlemma4} in order to show \eqref{yakirlem3}.
\begin{lemma}\label{lem1}
Let $\mu_X^\lambda$ be as in \eqref{muxlamb}. Then the diffraction measure $\widehat{\gamma_X^\lambda}=\tau_\lambda(\widehat{\gamma_{X}})$, where $\tau_\lambda$ is the action of the translation by $\lambda$; in particular, under (H.2) $\widehat{\gamma_X^\lambda}$ is singular with respect to Lebesgue measure.
\end{lemma}
We just sketch the proof, since the methods are standard.
\begin{proof}
Note that $f_\lambda(x)=e^{-2\pi i \langle x,\lambda\rangle}$ is bounded and $\mu_X$ is translation-bounded, therefore the measure $\mu_X^\lambda=f_\lambda \mu_X$ is translation-bounded as well, thus an autocorrelation measure $\gamma_X^\lambda$ exists (cf. \cite[Proposition 9.1]{BG1}). Now definition \eqref{weaklim} can be applied for $\mu=\mu_X^\lambda$, and by a change of variable in the formula \eqref{weaklim} for $\widehat{\gamma_\mu}$, the first claim follows. Finally, if $\widehat{\gamma_X}$ is singular with respect to the Lebesgue measure, the second claim follows directly since $\widehat{\gamma_{X}^{\lambda}}$ is just a translation of $\widehat{\gamma_X}$ by the vector $\lambda\in\mathbb{R}^d$.
\end{proof}
%%%%%%%%%%%%%%%%%%%%%%%%%%%%%%%%%%%%%%%%%%%%%%%%%%%%%%%%%%%%%%%%%%%%%%%%%%%%%%%%%%%%%%%%%%%%%%%%%%%%%%%%%%%%%%%%%%%%%%%%%%%%
\subsection{Independent random perturbations and proof of Theorem A}
\begin{prop}\label{prop2}
Let $X$ and $(\xi_p)_{p\in X}$ be as in Theorem A and let $\mu_\xi^\lambda$ be as in \eqref{muxilam}. Then almost surely, for all $\lambda\in \R^d$ there holds 
\begin{equation}\label{yakir1}
    \quad \gamma_\xi^\lambda = \mathsf{Var}\left[e^{-2\pi i \langle \xi, \lambda\rangle}\right]\mathsf{dens}(X)\delta_0.
\end{equation}
\end{prop}
\begin{proof}

Observe that $\gamma_\xi^\lambda$ is defined in duality with $C_c(\R^d)$ only, as follows:
\begin{equation*}
\forall f\in C_c(\R^d)\quad \int f(x)\ d\gamma_\xi^\lambda(x)= \lim_{R\to\infty}\frac{1}{Vol(B_R)}\sum_{p,q\in X\cap B_R}f(p-q) \mu_\xi^\lambda(p)\overline{\mu_\xi^\lambda(q)}.
\end{equation*}

We first prove the almost-sure convergence at fixed $\lambda$. For this we write
\begin{equation*}
    \sum_{p,q\in X\cap B_R}f(p-q) \mu_\xi^\lambda(p)\overline{\mu_\xi^\lambda(q)}=f(0)\sum_{p\in X\cap B_R} |\mu_\xi^\lambda(p)|^2 + \sum_{\substack{p,q\in X\cap B_R\\ p\neq q}}f(p-q) \mu_\xi^\lambda(p)\overline{\mu_\xi^\lambda(q)}.
\end{equation*}
Recall that in (H.1) we assumed that $\mathsf{dens}(X):=\lim_{R\to\infty}\frac{\#(X\cap B_R)}{Vol(B_R)}$ exists. Therefore, to prove \eqref{yakir1} it is sufficient to prove that almost surely for  $p_0$ in $X$ (whose precise choice do not matter, since the random variables $\xi_p$ are assumed to be i.i.d.) and for every compact set $K\subset\mathbb{R}^d$:
\begin{eqnarray}
\lim_{R\to\infty}\frac{1}{\#(X\cap B_R)}\sum_{p\in X\cap B_R}|\mu_\xi^\lambda(p)|^2 &=& \mathbb E[|\mu_\xi^\lambda(p_0)|^2],\label{1stclaim}\\[3mm]
\lim_{R\to\infty}\frac{1}{Vol(B_R)}\sum_{p\in X\cap B_R}\mu_\xi^\lambda(p)\sum_{\substack{q\in X\cap (K+p)\\ q\neq p}} \overline{\mu_\xi^\lambda(q)}&=&0.\label{2ndclaim}
\end{eqnarray}
By applying \eqref{2ndclaim} for $K$ a sub-level set for $\mathrm{Re}(f)_-$, $\mathrm{Re}(f)_+$, $\mathrm{Im}(f)_+$,  $\mathrm{Im}(f)_-$ and by the Cavalieri principle, it will imply that the sum against $f(p-q), p\neq q$ as well would tend to zero.

\medskip
{\bf Proof of \eqref{1stclaim}.} The random variables $\mu_\xi^\lambda(p),\mu_\xi^\lambda(q)$ are deterministic (bounded) functions of $\xi$, thus are also independent and identically distributed for $p\neq q$. We find \eqref{1stclaim} by a direct application of the classical strong law of large numbers.
\medskip

{\bf Proof of \eqref{2ndclaim}.} For $p_0,p_1\in X$ to be distinct fixed points, consider the $\mathbb P$-preserving maps $T_{p,q}:(\Omega, \mathbb P)\to(\Omega,\mathbb P)$ induced by the interchange of $p$ with $p_0$ and of $q$ with $p_1$, for $p\neq q\in X$ and 
\begin{equation}\label{tildear}
    \widetilde A_R\phi(\omega):=\frac{1}{Vol(B_R)}\sum_{\substack{p\in X\cap B_R\\ q\in X\cap(K+p)\setminus \{p\}}} \phi(T_{p,q}(\omega)).
\end{equation}
Our strategy of proof is to prove convergence of the $A_R\phi$ as $R\to \infty$, after which we can take $\phi(\omega):=\mu_\xi^\lambda(p_0)\overline{\mu_\xi^\lambda(p_1)}$ to conclude the first equality in \eqref{2ndclaim}.

\medskip
Now, we partition the averages from \eqref{tildear} over a collection of periodic subsets covering $\mathbb R^d$, as follows. For each positive integer $N$, for $r$ being the separation constant of $X$ we set
\[
C:=\left[-\frac{r}{2\sqrt d},\frac{r}{2\sqrt d}\right)^d,
\]

and for $i\in\{1\ldots,N-1\}^d$ we define the set $A_N^i$ by 

\[
A_N^i:=\bigcup_{k\in\frac{Nr}{\sqrt d}\mathbb Z^d}\left(C+ i\frac{r}{\sqrt d} + k\right).
\]

Note that the sets $\{A_N^i: i\in\{0,\dots,N-1\}^d\}$ form a partition of $\mathbb R^d$, and due to the fact that $r$ is the separation radius of $X$, for any $p\neq q\in A_N^i\cap X$ we have $p-q> (N-1)r/\sqrt d$. Assuming $N$ is such that $\mathrm{diam}(K)<(N-1)r/\sqrt d$ it follows $q\notin (K+p)\setminus \{p\}$. Therefore, for such $N$ the random variables
\begin{equation*}\label{rvain}
    \left(\sum_{\substack{q\in X\cap(K+p)\setminus \{p\}}} \mu_{\xi}^{\lambda}(p)\overline{\mu_{\xi}^{\lambda}(q)}\right)_{p\in X\cap A_N^i}
\end{equation*}
are independent. Furthermore, we claim that under the above hypotheses on $K,N$, if $\phi$ is bounded and $\mathbb E[\phi]=0$ then the following holds: if we define
\begin{equation*}
A_{R,N}^i\phi(\omega):=\frac{1}{Vol(B_R)}\sum_{\substack{p\in X\cap B_R\cap A_N^i\\ q\in X\cap(K+p)\setminus \{p\}}} \phi(T_{p,q}(\omega)),
\end{equation*}

then $A_{R,N}^i\phi(\omega)$ almost surely converges to $0$, i.e

\begin{equation}\label{convairn}
\mathbb P\left\{\omega:\ \lim_{R\to\infty}A_{R,N}^i\phi(\omega)=0\right\}=1 . 
\end{equation}

Indeed, we can take $\phi_j$ to be independent copies of $\phi$ and write equivalently
\begin{equation*}
A_{R,N}^i\phi=\frac{n_{R,N}^i}{Vol(B_R)}\cdot\frac1{n_{R,N}^i}\sum_{j=1}^{n_{R,N}^i} \phi_j, \qquad n_{R,N}^i:=\#\left\{(p,q)\in X^2:\ \left.\begin{array}{l}p\in B_R\cap A_N^i\\ q-p\in K\setminus \{0\}\end{array}\right.\right\}.
\end{equation*}
First note that since $X$ is uniformly discrete and $K$ is bounded, we have that 
\begin{equation}\label{bdnpts}
\sup_{R>0}\frac{n_{R,N}^i}{Vol(B_R)}<\infty.
\end{equation}
Then for each subsequence $R_k\to \infty$ such that $n_{R_k,N}^i/Vol(B_{R_k})$ converges there are two possibilities: 
\begin{enumerate}
    \item The sequence $n_{R_k,N}^i/ Vol(B_{R_k})\to 0$, in which case by triangle inequality and using the fact that $\phi$ is bounded, we find \eqref{convairn}.
    \item We have %$0<\liminf_{k\to\infty}n_{R_k,N}^i/Vol(B_{R_k})$,
    $0<\lim_{k\to\infty}n_{R_k,N}^i/Vol(B_{R_k})$, in particular $n_{R_k,N}^i\to \infty$. By boundedness of the sequence $(\phi_j)_{j}$ and by the strong law of large numbers for non-necessarily identically distributed random variables (see for instance \cite[Corollary 1]{Ch-Gos}), we get  $\frac{1}{n_{R_k,N}^i}\sum_{j=1}^{n_{R_k,N}^i}\phi_j\to 0$. Using \eqref{bdnpts}, we obtain \eqref{convairn}
    also in this case.
\end{enumerate}
We note that $\mu_\xi^\lambda$ is bounded and that
\[
\mathbb E[\mu_\xi^\lambda(p_0)\overline{\mu_\xi^\lambda(p_1)}]=\mathbb E[\mu_\xi^\lambda(p_0)]\overline{\mathbb E[\mu_\xi^\lambda(p_1)]}=0\quad\text{ for }p_0\neq p_1\in X,
\]
using the fact that the random variables $\mu_\xi^\lambda(p_0), \mu_\xi^\lambda(q_0)$ are i.i.d. and of zero mean. We thus have \eqref{convairn} and by summing over $i$ we obtain the desired limit \eqref{2ndclaim}. 

\medskip

In order to show that the almost sure convergence holds also contemporarily for all $\lambda\in\mathbb R^d$, we proceed like in \cite{Yakir}: Firstly, we have that almost-sure convergence can be obtained at for a countable dense subset of $\lambda\in\R^d$. To extend the almost sure convergence to all $\lambda$, we prove the almost-sure (uniform in $\lambda$) control of the $\lambda$-gradients
\begin{equation}\label{lderbd}
\limsup_{R>0}\sup_{\lambda\in\mathbb R^d} \left|\frac1{Vol(B_R)}\sum_{p\in X\cap B_R}\sum_{q\in X\cap (K+p)}\nabla \left(\mu_\xi^\lambda(p)\overline{\mu_\xi^\lambda(q)}\right)\right|.
\end{equation}
For this we just apply triangular inequality, via the rough bounds
\[
    \left|\nabla \mu_\xi^\lambda(p)\right|\le \left|\xi_p\right| + \mathbb{E}[|\xi|], \qquad \left|\mu_\xi^\lambda(p)\right|\le 2.
\]
We thus find for $\phi:=\mu_\xi^\lambda(p_0)\overline{\mu_\xi^\lambda(q_0)}$:
\[
\limsup_{R\to\infty}\sup_{\lambda\in\mathbb R^d}\left|\nabla \widetilde A_R\phi(\omega)\right|\le \limsup_{R\to\infty}\frac2{Vol(B_R)}\sum_{\substack{p\in X\cap B_R\\ q\in X\cap (K+p)\setminus \{p\}}}(|\xi_p(\omega)|+|\xi_q(\omega)|+2\mathbb E[|\xi|]).
\]

By partitioning $K$ if it is necessary, we can consider $K$ of small diameter (which does not change the conclusion of \eqref{2ndclaim}), and then $X\cap (K+p)$ has cardinality at most $1$ for all $p$. Thus 
\begin{equation}\label{lastbdgradient}
\sum_{\substack{p\in X\cap B_R\cap A_N^i\\ q\in X\cap (K+p)\setminus \{p\}}}(|\xi_p(\omega)|+|\xi_q(\omega)|+2\mathbb E[|\xi|])\le 2\sum_{p\in X\cap (B_R+K)}(|\xi_p|+\mathbb E[|\xi|]),
\end{equation}
and now we can label points $\{p_n\}_{n\in\N}= X\cap A_N^i$ so that $|p_n|$ is increasing, and apply the classical strong law of large numbers to the random variables $|\xi_{p_n}|$ to show that the %above
right hand side in \eqref{lastbdgradient} is almost surely bounded by $C_X\mathbb E[|\xi|]Vol(B_R)$. This gives that almost surely \eqref{lderbd} is bounded by $C_X\mathbb E[|\xi|]$ in which $C_X$ is a packing constant, $C_X\le \frac{C_d}{r^d}$ in which $C_d$ depends only on the dimension and $r>0$ is the uniform separation radius of $X$.

\medskip

Note that, if we fix a countable dense $\Lambda\subset\mathbb{R}^d$, almost surely \eqref{1stclaim} and \eqref{2ndclaim} hold for every $\lambda\in\Lambda$. Thus, to get the almost sure convergence in \eqref{2ndclaim} contemporarily also for all $\lambda\in\R^d\setminus\Lambda$, observe that if $(\lambda_n)_{n\in\mathbb{N}}\subset\Lambda$ is a sequence such that $\lambda_n\longrightarrow\lambda$ , then for  $\phi^{\lambda}=\mu_{\xi}^{\lambda}(p_0)\overline{\mu_{\xi}^{\lambda}(q_0)}$ there holds

\begin{equation*}
\begin{split}
    \left|\widetilde{A}_R\phi^{\lambda}\omega\right|&\leq \left|\widetilde{A}_R\phi^{\lambda}\omega-\widetilde{A}_R\phi^{\lambda_n}\omega\right|+\left|\widetilde{A}_R\phi^{\lambda_n}\omega\right|\\
    &\leq M_K|\lambda-\lambda_n|+\left|\widetilde{A}_R\phi^{\lambda_n}\omega\right|,
\end{split}
\end{equation*}

where $M_K$ is the value of \eqref{lderbd}. By taking the limsup as $R\to\infty$ and then $n\to\infty$ in the last inequality, we obtain that on the same event on which convergence holds for the $\lambda_n$ we also have $\widetilde{A}_R\phi^{\lambda}\longrightarrow 0$. Thus this convergence holds almost surely for every $\lambda\in\mathbb{R}^d$. Finally, an analogous argument together an application of the $\lambda$-continuity of $\mu_{\xi}^{\lambda}$ and the dominated convergence theorem, we get the convergence in \eqref{1stclaim} almost surely for all $\lambda$. This concludes the proof of Proposition \ref{prop2}.
\end{proof}

\begin{proof}[Proof of Proposition \ref{reduction Main Thm}.]
From Theorem \ref{yakirlemma4}, Lemma \ref{lem1} and Proposition \ref{prop2} (resp. Propositions \ref{noind-prop2} and \ref{stationarycase} below in the case of Theorems B and C), we have that for every positive test function $f\in\mathcal{S}(\mathbb{R}^d)$

\begin{equation}\label{limzero}
    \lim_{R\to\infty}\dfrac{1}{Vol(B_R)}\left|\displaystyle\sum_{p\in X\cap B_R}\mu_{\xi}^{\lambda}(p)\overline{\mu_X^{\lambda}(p)}\widehat{f}(0)+\sum_{\substack{p,q\in X\cap B_R\\ q\neq p}}\mu_{\xi}^{\lambda}(p)\overline{\mu_X^{\lambda}(q)}\widehat{f}(p-q)\right|=0.
\end{equation}

In particular, we can choose $f=f_\epsilon\in\mathcal{S}(\mathbb{R}^d)$ such that $1<\widehat{f_{\varepsilon}}(0)$ (in effect, what we use below is that $0<\widehat{f_\varepsilon}(0)$ and $\varepsilon/\widehat{f_\varepsilon}(0)\to 0$, for $\epsilon\to 0$) and

\begin{equation*}
    \forall p\in X\qquad\sum_{q\in X\setminus\{p\}}|\widehat{f}_\epsilon(p-q)|<\epsilon.
\end{equation*}

This is possible by the uniform discreteness of $X$, by taking $\hat f_\epsilon$ radially decreasing and concentrated on $B(0, r/2)$ where $r$ is the separation constant of $X$ in the sense that $\int_{\mathbb R^d\setminus B(0,r/2)} |\hat f_\epsilon|< \frac{r^d\epsilon}{C}$. For example, this holds for $\hat f_\epsilon$ a Gaussian centered at $0$ with small enough variance. For such $f_\epsilon$ we can estimate $\sum_{q\in X\setminus \{p\}}|\hat f_\epsilon(p-q)|$ by the abovementioned integral, using the separation of $X$ and the fact that $\hat f_\epsilon$ is radially decreasing.

Now fix $R_0$ such that the quantity in \eqref{limzero}, is smaller than $\epsilon$ and that $\frac{\#(X\cap B_R)}{Vol(B_R)}<\mathrm{dens}(X)+1$ for all $R\geq R_0$. By triangle inequality applied to \eqref{limzero}, for such $R$ we get

\begin{equation*}
    \dfrac{\widehat{f_{\varepsilon}(0)}}{Vol(B_R)}\left|\displaystyle\sum_{p\in X\cap B_R}\mu_{\xi}^{\lambda}(p)\overline{\mu_X^{\lambda}(p)}\right|\leq \varepsilon+\dfrac{2}{Vol(B_R)}\sum_{\substack{p,q\in X\cap B_R\\ q\neq p}}|\widehat{f}_\epsilon(p-q)|\le (2\ \mathrm{dens}(X)+3)\ \epsilon,
\end{equation*}

which follows by noting that $\sup_{p\neq q\in X}|\mu_\xi^\lambda(p)\overline{\mu_X^\lambda(q)}|\le 2$. By dividing by $\widehat{f_{\varepsilon}(0)}$, we then take $\epsilon\to 0$ in the above inequality, and this completes the proof of Proposition \ref{reduction Main Thm} in the case of i.i.d random perturbations.
\end{proof}

%%%%%%%%%%%%%%%%%%%%%%%%%%%%%%%%%%%%%%%%%%%%%%%%%%%%%%%%%%%%%%%%%%%%%%%%%%%%%%%%%%%%%%%%%%%%%%%%%%%%%%%%%%%%%%%%%%%%%%%%%%%%%%%%%%%%%%%%%%%%%%%%%%%%%%%%%%%%%%%%%%%%%%%%%%%%%%%%%%%%%%%%%%%%%%%%%%%%%

\subsection{Non-independent random perturbations and proof of Theorem B}\label{ssec: non-ind rv}
For extending Proposition \ref{prop2} to the case of non-independent identically distributed $\xi_p,p\in X$, we need a mixing hypothesis. Several notions of mixing are available (see the survey \cite{bradleymix}); we perform our proofs based on the following ad-hoc strong mixing condition. For $p\in X$ and $K\subset \R^d$ compact, define the sub-$\sigma$-algebras 
\[
\mathcal F_p:=\sigma(\xi_p),\quad \mathcal F_{p,K}:=\sigma(\xi_q:\ q \in (K+p)\cap X\setminus \{p\}).
\]
we say that the random process $(\xi_p)_{p\in X}$ is {\em compactly-mixing} if for any fixed compact $K\subset \R^d$ there holds
\begin{equation}\label{mix-salg}
   \lim_{R\to\infty}\frac1{Vol(B_R)}\sum_{p\in X\cap B_R} \sup\left|\mathbb {E}[fg]- \mathbb {E}[f]\mathbb {E}[g]\right|=0.
\end{equation}
where for each $p\in X$ the supremum is taken over the set of all functions $f\in L^{\infty}(\mathcal{F}_p,\mathbb{P})$, $g\in L^{\infty}(\mathcal{F}_{p,K},\mathbb{P})$ such that $||f||_{L^{\infty}}\leq 1$, $||g||_{L^{\infty}}\leq 1$. Observe that \eqref{mix-salg} is equivalent to 

\begin{equation*}
    \lim_{R\to\infty}\frac1{Vol(B_R)}\sum_{p\in X\cap B_R} \sup\left|\mathbb {P}[A\cap B]- \mathbb {P}[A]\mathbb{P}[B]\right|=0,
\end{equation*}

where for each $p\in X$, the supremum is taken over the events $A\in \mathcal{F}_p$ and $B\in\mathcal{F}_{p,K}$.\\

We first prove the following:
\begin{lemma}
For $\mu_{\xi}^{\lambda}$ as before and for $\lambda\in\mathbb{R}^d$ fixed, we have that almost surely
 \begin{equation}\label{eqslln}
    \frac{1}{Vol(B_R)}\sum_{p\in X\cap B_R}\sum_{\substack{q\in X\cap (K+p)\\ q\neq p}}\left[\mu_{\xi}^{\lambda}(p)\overline{\mu_{\xi}^{\lambda}(q)}-\mathbb{E}[\mu_{\xi}^{\lambda}(p)\overline{\mu_{\xi}^{\lambda}(q)}]\right]\longrightarrow 0.
\end{equation}
\end{lemma}

\begin{proof}
To prove this, we apply Proposition \ref{prop:slln} to the real and imaginary parts of the variables $\mu_\xi^\lambda(p)\overline{\mu_\xi^\lambda(q)}$ corresponding to elements $p,q\in X$ such that $p-q\in K$. Firstly, by (H.1) we may replace the factor $1/Vol(B_R)$ by $1/\#(X\cap B_R)$ in \eqref{eqslln}. Note that the sets of index pairs figuring in the double sums in \eqref{eqslln} form a sequence increasing in $R$, thus we can find a labelling of $X_i:=\mu_\xi^\lambda(p)\overline{\mu_\xi^\lambda(q)}$ so that the convergence in \eqref{eqslln} follows from the convergence of a subsequence of the $\frac1n\sum_{n=1}^n(X_k -\mathbb E[X_k])$. We apply Proposition \ref{prop:slln} to obtain convergence of the full sequence. We verify the hypotheses of Proposition \ref{prop:slln} next, for this case:

\begin{itemize}
    \item The bounds $C_k$ and the control of part 1 of the hypotheses of Proposition {\ref{chgos}} follow because $|\mu_\xi^\lambda(p)|\le 2$ almost surely.
    
    \item  For the summability of variances along geometric series of hypothesis 2 of Proposition \ref{chgos}, it suffices to prove $\frac{\mathsf{Var}(S_{k_n})}{k_n}\le C$, which follows from the boundedness of $\mathsf{Var}(\mu_{\xi}^{\lambda}(p)\overline{\mu_{\xi}^{\lambda}(q)})$ and $\mathsf{Cov}(\mu_{\xi}^{\lambda}(p),\mu_{\xi}^{\lambda}(q))$. In fact, for $X_i=\mu_{\xi}^{\lambda}(p)\overline{\mu_{\xi}^{\lambda}(q)}$ as before, we get:
\begin{equation*}
    \mathsf{Var}(S_{k_n})=\sum_{i=1}^{k_n}\mathsf{Var}(X_i)+2\sum_{\substack{i,j=1\\ i\neq j}}^{k_n}\mathsf{Cov}(X_i,X_j), 
\end{equation*}
and by the definition of $\mu_{\xi}^\lambda(p)$ we have
\begin{equation*}
\begin{split}
\mathsf{Var}(X_i)&\leq\mathbb{E}\left[\left|\mu_{\xi}^{\lambda}(p)\overline{\mu_{\xi}^{\lambda}(q)}\right|^2\right]\leq 16\\
|\mathsf{Cov}(X_i,X_j)|&=\left|\mathbb{E}\left(\mu_{\xi}^{\lambda}(p)\overline{\mu_{\xi}^{\lambda}(q)}\overline{\mu_{\xi}^{\lambda}(p')}\mu_{\xi}^{\lambda}(q')\right)-\mathbb{E}\left[\mu_{\xi}^{\lambda}(p)\overline{\mu_{\xi}^{\lambda}(q)}\right]\mathbb{E}\left[\overline{\mu_{\xi}^{\lambda}(p')}\mu_{\xi}^{\lambda}(q')\right]\right|\leq 32.
\end{split}
\end{equation*}
\end{itemize}

Therefore, the limit (\ref{eqslln}) follows as a consequence of Proposition \ref{prop:slln} as claimed.
\end{proof}

\begin{lemma}
If the random process $(\xi_p)_{p\in X}$ is compactly-mixing, then there holds
\begin{equation}\label{discrepancycontrol}
    \sum_{p\in X\cap B_R}\sum_{\substack{q\in X\cap (K+p)\\ q\neq p}}\mathbb{E}[\mu_{\xi}^{\lambda}(p)\overline{\mu_{\xi}^{\lambda}(q)}]=o(R^d).
\end{equation}
\end{lemma}
\begin{proof}
Indeed, this a direct consequence of the definition \eqref{mix-salg}, observing that $\mu_\xi^\lambda(p)$ is almost surely bounded for all $\lambda$.
\end{proof}
\medskip
Now \eqref{eqslln} and \eqref{discrepancycontrol} give the following:
\begin{prop}\label{noind-prop2}
Assume that the $(\xi_p)_{p\in X}$ are identically distributed and compactly-mixing, and let $\mu_\xi^\lambda$ be as in \eqref{muxilam}. Then almost surely, for all $\lambda\in\mathbb R^d$ there holds

\begin{equation*}
    \gamma_{\xi}^{\lambda}=\mathsf{Var}\left[e^{-2\pi i\langle\xi,\lambda\rangle}\right]\mathsf{dens}(X)\delta_0.
\end{equation*}
\end{prop}
To obtain Proposition \ref{reduction Main Thm} under the above hypotheses, we now note that the upper bound of \eqref{lderbd} can also be extended under our assumed bound on the covariances of the $\xi_p$. Indeed all the proof up to \eqref{lastbdgradient} can be repeated verbatim, and for the bound of the last term in \eqref{lastbdgradient}, we may now apply Proposition \ref{Etem-Cov} to the $|\xi_p|$ after to enumerate the points in $X\cap B_R$ and by using (H.1), which gives a law of strong numbers under our assumed condition \eqref{covhyp} on the covariances of the $|\xi_p|$. This concludes the proof of Theorem B. 

\medskip
\begin{example}[{\em A recipe to produce compactly-mixing random processes in $X$}]
Let $(R_k)_{k\geq 1}$ be a sequence of positive real numbers such that $R_k,R_k/R_{k-1}\longrightarrow\infty$, and denote $\Omega_k:=B_{R_k}\setminus B_{R_{k-1}}$ for $k\ge 2$. Consider a process of identically distributed random vectors $(\xi_p)_{p\in X}$ such that the following hypotheses on correlations hold.
\begin{enumerate}
    \item For every $p\in X\cap B_{R_1}$, the random vector $\xi_p$ has no correlation with the $\xi_{p'}, p'\in X\setminus\{p\}$.
    \item For every even $k>1$, the $\xi_p, p\in X\cap \Omega_k$ are not correlated to the remaining $\xi_{p'}, p'\in X\setminus \Omega_k$.
    \item For every odd $k>1$ and for for $p\neq p'\in X\cap \Omega_k$, $\xi_p$ and $\xi_{p'}$ are uncorrelated, and for every $p\in X\cap \Omega_k$, there exists at least one $q\in X\cap \Omega_{k-1}$ such that $\xi_p $ has nontrivial correlation to the value of $\xi_q$. 
\end{enumerate}
It is not hard to check that this procedure produces a random process satisfying (\ref{mix-salg}), since correlations at distance closer than $\mathrm{diam}(K)$ do not occur. %
\end{example}

As a further remark, which allows to understand the role of the hypothesis that the $\xi_p$ are compactly-mixing, which is used here in order to obtain \eqref{discrepancycontrol}. The following result makes it explicit, that the outcome of this is to obtain a version of \eqref{yakir1}, but in expectation rather than almost-surely. As a consequence, we could also require as a hypothesis the control stated \eqref{yakir1exp} below, rather than the compact-mixing hypothesis.
\begin{lemma}\label{expect-yak}
If $\xi_p$ are identically distributed and if $X$ verifies (H1) then the fact that \eqref{discrepancycontrol} holds for all $\lambda\in \mathbb R^d$ is equivalent to the following:
\begin{equation}\label{yakir1exp}
\forall \lambda\in \mathbb R^d,\qquad \mathbb E(\gamma_\xi^\lambda)=\mathsf{Var}\left[e^{-2\pi i\langle \xi,\lambda\rangle}\right]\mathrm{dens}(X)\delta_0.
\end{equation}
\end{lemma}
\begin{proof}
Indeed, we have
\[
\lim_{R\to\infty} \frac1{Vol(B_R)}\sum_{p\in X\cap B_R} \mathbb E\left[|\mu_\xi^\lambda(p)|^2\right]=\mathrm{dens}(X)\mathsf{ Var}(\mu_\xi^\lambda),
\]
thus \eqref{discrepancycontrol} is equivalent to
\begin{equation}\label{dc2}
    \sum_{p\in X\cap B_R}\sum_{q\in X\cap (K+p)}\mathbb{E}\left[\mu_{\xi}^{\lambda}(p)\overline{\mu_{\xi}^{\lambda}(q)}\right]=\mathrm{dens}(X)\ Vol(B_R)\mathsf{Var}(\mu_\xi^\lambda)+o(R^d).
\end{equation}
The fact that \eqref{dc2} is valid for all compact $K$ is equivalent to the convergence of the corresponding $(Vol(B_R))^{-1}$-rescaled sum tested against all test functions $f$, which is the weak formulation of \eqref{yakir1exp}.
\end{proof}

\subsection{Extension of Theorems A, B for $X$ being the asymptotically affine deformation of a lattice}\label{asyaffine}
In this section we extend Theorems A, B by weakening the independence requirements on the random perturbations $\xi_p$. To achieve this we limit ourselves to a smaller class of sets $X$. 

\medskip
A useful tool to replace the independence of the $\xi_p$, $p\in X$, is the notion of {\em ergodicity} of the process $(\xi_p)_{p\in X}$, labelled by the set $X$. To make sense of this notion, commonly $X$ is required to have more structure, e.g a transitive group action on $X$ is required, and classically $X$ itself is assumed to be a lattice. Therefore, firstly we study the case of $X$ being a lattice: in this case we recover the main results in \cite{Yakir}. Secondly, we consider the case when $X$ is almost-affine, meaning that though $X$ lacks a group structure (the sum of elements of $X$ is not generally in $X$), it is a small perturbation of a subgroup of $\R^d$. In this case the results still depend strongly on a group structure. For a discussion of further weakenings on the structure of $X$, see Section \ref{stationarygen}.

\subsubsection{Lattice case}\label{sec-latcase}
As was pointed out by Yakir in \cite{Yakir}, in the case of $X$ being a lattice, the condition that the process $(\xi_p)_{p\in X}$ is ergodic is sufficient for (\ref{autocorr}) to converge as $R\to\infty$.
\medskip

\begin{lemma}[Lattice case]\label{latcase}
If $X$ is a lattice and $(\xi_p)_{p\in X}$ is stationary and ergodicwith respect to lattice shifts, then almost surely
\begin{equation}\label{gammalat}
    \gamma_{\xi}^{\lambda}=\sum_{k\in X}\mathbb E\left[\mu_\xi^\lambda(0)\overline{\mu_{\xi}^{\lambda}(k)}\right]\mathrm{dens}(X)\delta_k.
\end{equation}
\end{lemma}
\begin{proof}
By stationarity and using Birkhoff's theorem for the natural shifts on $X$, we ensure the existence of the limit defining $\gamma_\xi^\lambda(k)$ as in \eqref{autocorr}. Ergodicity of the process $(\xi_p)_{p\in X}$ yields that the above limit is almost surely constant and thus almost surely converges to $\mathbb{E}[\mu_{\xi}^{\lambda}(0)\overline{\mu_{\xi}^{\lambda}(k)}]$. 
\end{proof}
In general a measure like $\gamma_\xi^\lambda$ as given in \eqref{gammalat} may have Fourier transform $\widehat{\gamma_\xi^\lambda}$ with nonzero atomic part (i.e. the initial $\sum_{p\in X}\mu_\xi^\lambda(p)\delta_p$ may have pure-point spectrum). If however we strengthen the ergodicity condition on the map $(\xi_p)_{p\in X}$ to a \emph{weak mixing} condition, this is precluded. We require the following weak mixing condition. Let $\tau_k: X\to X$ be the translation by $k\in X$. Associated to each $\tau_k$ we define a measurable $\mathbb P$-preserving transformation $T_k$, which is induced by interchanging the $X$-coordinates $(\xi_p)_{p\in X}$ according to $\tau_k$. The process $(\xi_p)_{p\in X}$ is weakly mixing if for all $f,g\in L^2(\mathbb{P})$ there holds

\begin{equation}\label{wmixL2}
\lim_{R\to\infty}\frac{1}{Vol(B_R)}\sum_{k\in X\cap B_R}\left|\mathbb E[f\circ T_k\ \overline g] - \mathbb E[f]\overline{\mathbb E[g]}\right|=0.
\end{equation}

Note that the compactly-mixing condition introduced earlier in \eqref{mix-salg} is not directly comparable to \eqref{wmixL2}, and we do not pursue an in-depth comparison in this work.
\medskip

We have the following
\begin{lemma}[Basic lattice case]\label{bala_case}
If $X$ is a lattice and $(\xi_p)_{p\in X}$ is stationary and weakly mixing then almost surely for every $\lambda\in\R^d$ the random diffraction measure $\widehat{\gamma_\xi^\lambda}$ has no atoms.
\end{lemma}

\begin{proof}
Recall Strungaru's criterion \cite[Proposition 4.1]{Strungaru} (here formulated in the setting of formula \eqref{gammalat} which we know to hold for $\gamma_\xi^\lambda$): 
\begin{equation}\label{strungaru}
\lim_{R\to\infty}\frac1{Vol(B_R)}\sum_{k\in X\cap B_R}\left|\mathbb E[\mu_\xi^\lambda(0)\overline{\mu_\xi^\lambda(k)}]\right| =0\quad\Leftrightarrow\quad \widehat{\gamma_\xi^\lambda} \text{ has no atoms}.
\end{equation}
We note that \eqref{strungaru} follows from \eqref{wmixL2}, if we apply it for $f(\omega)=g(\omega)=e^{2\pi i\langle \lambda, \xi_0(\omega)\rangle}$, since by stationarity $\xi_k=\xi_0\circ T_k$ and
\begin{equation*}\label{corr1}
\mathbb E\left[\mu_\xi^\lambda(0)\overline{\mu_\xi^\lambda(k)}\right]=\mathbb E\left[e^{-2\pi i\langle \xi_0-\xi_k, \lambda\rangle}\right] - \mathbb E\left[e^{-2\pi i\langle \xi_0,\lambda\rangle}\right]\overline{\mathbb E\left[e^{-2\pi i\langle \xi_k,\lambda\rangle}\right]}.
\end{equation*}
\end{proof}

\subsubsection{Asymptotically affine perturbations of a lattice}
Now assume that $X=\psi(L)$ where $L$ is a lattice, and that $\xi_p$ are identically distributed and satisfy $\mathbb E(|\xi_p|^{d+\epsilon})<\infty$ as before. Denote $\mu'(x):=\mu_\xi^\lambda(\psi(x))$ for $x\in L$, and assume that $\mu'(x)$ are stationary and ergodic under the natural shifts of $x\in L$. We have the following asymptotic equalities, valid in the limit as $R\to \infty$, with error estimates as in Lemma \ref{controlled mass}:
\begin{eqnarray}
    \lefteqn{\frac1{Vol(B_R)}\sum_{p,q\in X\cap B_R} g(p-q)\mu_\xi^\lambda(p)\overline{\mu_\xi^\lambda(q)}}\nonumber\\
    &\simeq&\frac1{Vol(B_R)}\sum_{x,y\in L\cap B_R}g(\psi(x)-\psi(y))\mu'(x)\overline{\mu'(y)}\nonumber\\
    &\simeq&\sum_{k\in L\cap \mathrm{supp}(g)}\frac1{Vol(B_R)}\sum_{x\in L\cap B_R} g(\psi(x)-\psi(x+k)) \mu'(x)\overline{\mu'(x+k)}%\nonumber\\
    \label{bdcaseq}
\end{eqnarray}
We find it natural based on \eqref{bdcaseq} to introduce the following notion. We say that $\psi$ has {\it asymptotically affine} if  one of the following equivalent conditions holds:
\begin{enumerate}
\item There exists $F:L\to \mathbb R$ such that for any continuous compactly supported $g$ and for any $k\in L$
\begin{equation}\label{aus}
\lim_{R\to \infty} \frac1{Vol(B_R)}\sum_{x\in L\cap B_R}g(\psi(x)-\psi(x+k))= g(F(k)).
\end{equation}
\item For all $k\in L$ there holds
\begin{equation}\label{asyaff}
\lim_{\substack{x\in L\\|x|\to\infty}}(\psi(x)-\psi(x+k))=F(k).
\end{equation}
\item The limit \eqref{asyaff} holds for a linear function $F$.
\end{enumerate}

To see that $1\Rightarrow 2$ above, it suffices to test \eqref{aus} over a sequence of $g_n$ that approximate a Dirac mass at $F(k)$. Then $2\Rightarrow 3$ follows by observing that if the limit in \eqref{asyaff} exists for $k\in\{k_1,k_2\}$ then \[F(k_1)+F(k_2)=\lim_{\substack{x\in L\\|x|\to\infty}}(\psi(x)-\psi(x+k_1))+\lim_{\substack{x\in L\\|x|\to\infty}}(\psi(x+k_1)-\psi(x+k_1+k_2))= F(k_1+k_2),
\]
and finally, $3\Rightarrow 1$ is a direct verification.

\medskip

Now if we assume that $\psi$ satisfies \eqref{aus} and we insert this limit into \eqref{bdcaseq}, using the fact that $\mu'(x)$ are bounded and that $\mathrm{supp}(g)\cap L$ is finite, we directly find
\begin{eqnarray}\label{bdcaseq2}
\text{\eqref{bdcaseq}}&\simeq&\sum_{k\in L\cap \mathrm{supp}(g)}\frac1{Vol(B_R)}\sum_{x\in L\cap B_R} g(F(k)) \mu'(x)\overline{\mu'(x+k)}\nonumber \\
&\simeq&\sum_{k\in L\cap \mathrm{supp}(g)}\mathbb E[\mu_\xi^\lambda(0)\overline{\mu_\xi^\lambda(k)}]\mathrm{dens}(L)g(F(k)),
\end{eqnarray}
where in the second step we applied Lemma \ref{latcase}.
We thus have from \eqref{bdcaseq} directly the following result, which extends Lemmas \ref{latcase} and \ref{bala_case}:

\begin{prop}[asymptotically affine perturbations of a lattice]\label{bdlat1}
Assume that $X=\psi(L)$ where $L$ is a lattice and $\psi$ is asymptotically affine with limit linear part $F$. Assume that the random process $\xi_{\psi(x)}$, $x\in L$ is stationary and ergodic under shifts of $L$, and that the $\xi_{\psi(x)}$ are identically distributed. Then the limit defining $\gamma_\xi^\lambda$ as in \eqref{autocorr} exists and we have
\begin{equation}\label{latcase2}
\gamma_\xi^\lambda=\sum_{k\in L}\mathbb E\left[\mu_\xi^\lambda(0)\overline{\mu_\xi^\lambda(k)}\right]\mathrm{dens}(L)\delta_{F(k)}.
\end{equation}
Furthermore, if $\xi_{\psi(x)}, x\in L$ is weakly-mixing under shifts of $L$, then $\widehat{\gamma_\xi^\lambda}$ has no atoms.
\end{prop}
The above result directly allows to obtain Proposition \ref{reduction Main Thm}, which as before, together with Lemma \ref{controlled mass} gives the following result:
\medskip

{\bf Theorem D. }\label{Thm_d}{\em For $X=\psi(L)\subset \R^d$ as in Proposition \ref{bdlat1}, assume that the random process $\xi_{\psi(x)}$, $x\in L$ is stationary and weakly mixing under shifts of $L$, that the $\xi_p$'s are identically distributed and that there exists a positive number $\varepsilon$ such that $\mathbb{E}(|\xi|^{d+\varepsilon})<\infty$. Then, almost surely for every $\lambda\in\mathbb{R}^d$, in the sense of vague convergence of measures, the same limit \eqref{eqthmA} as in Theorem A holds.}

\subsection{Stationary random perturbations and proof of Theorem C}

\subsubsection{Stationary random measures}
The proof of Theorem C combines the recovery problem setting in which $X_\xi=\{p+\xi(p):\ p\in X\}$ as before, with the study of \cite{gouere} who established a diffraction theory for stationary point processes. The setting of point processes is not sufficient for our setting, therefore we first generalize Gou\'er\'e's work to locally finite random measures.  

\medskip

We start with a series of definitions and properties that are classical for point processes, see e.g. \cite{neveu}. The case of point processes is the special case in which the random measures $\mu$ below are the empirical measures of random sets $X\in \R^d$, i.e. $\mu=\sum_{x\in X}\delta_x$. Proofs valid for this special case directly generalize to our setting and are well-known, therefore we only sketch the proofs.
 
\medskip

For any $f\in C_c(\mathbb{R}^d)$ define the functions $N_f:\mathcal{M}^{\infty}\to\mathbb{C}$ by $N_f(\mu):=\mu(f)=\langle\mu,f\rangle$. Let $\mathcal{A}$ be the $\sigma$-algebra on $\mathcal{M}^{\infty}$ generated by the $N_f$'s. 
A {\em random measure} is a measurable function $\mu:\Omega\to\mathcal{M}^{\infty}$ defined on a probability space $(\Omega,\mathcal{F},\mathbb{P})$ with values in $(\mathcal{M}^{\infty},\mathcal{A})$. We say that a random measure is {\em integrable} if the functions $N_f$ are integrable. The random measure is called {\em stationary} if for every $x_1,\ldots,x_k\in\mathbb{R}^d$ and every $F_1,\ldots,F_k\in \mathcal A$,
\begin{equation*}
    \mathbb{P}\{\omega\in\Omega:\ \mu*\delta_{x_1+t}\in F_1,\ldots,\mu*\delta_{x_k+t}\in F_k\}=\mathbb{P}\{\omega\in\Omega:\ \mu*\delta_{x_1}\in F_1,\ldots,\mu*\delta_{x_k}\in F_k\}
\end{equation*}
holds for every $t\in\mathbb{R}^d$.\\

Let $\mu$ be a stationary and integrable random measure. The {\em Palm measure of $\mu$} is the measure $\widetilde{P}$ on $(\mathcal{M}^{\infty},\mathcal{A})$ defined by

\begin{equation*}
    \widetilde{P}(F):=\frac{1}{Vol(B_1)}\cdot\mathbb{E}\left(\sum_{x\in \mathsf{supp}(\mu)\cap B_1}\mathbbm{1}_F(\mu*\delta_{-x})\right),\hspace{0.5cm} F\in\mathcal{A}.
\end{equation*}

The {\em intensity} $I(\widetilde{P})$ of the Palm measure $\widetilde{P}$ is the functional on $C_c(\mathbb{R}^d)$ defined by

\begin{equation*}
    I(\widetilde{P})(f):=\int_{\mathcal{M}^{\infty}}\mu(f)d\widetilde{P}(\mu).
\end{equation*}

\begin{lemma}[{Cf. \cite[Prop. 2.24]{neveu}}]
 Let $\mu$ be a stationary and integrable random measure. The intensity $I(\widetilde{P})$ of its Palm measure is a translation-bounded (and therefore tempered) measure.
\end{lemma}

Let $\psi\in C_c(\mathbb{R}^{d})$ be a non-negative function such that $\mathsf{supp}(\psi)=B_M$, $\int\psi(x)dx=1$, and which will be kept fixed from now on. 
Let $\mathcal{M}^{\infty}_{X}$ be the set of measures $\mu\in\mathcal{M}^{\infty}$ for which there is $t\in\mathbb{R}^d$ such that $\mathsf{supp}(\mu)= X+t$.
For every $f\in C_c(\mathbb{R}^d)$, define the function $H_f:\mathcal{M}^{\infty}_{X}\to\mathbb{C}$ by

\begin{equation*}
    H_f(\nu):=\sum_{x,y\in\mathsf{supp}(\nu)}\psi(x)\nu(x)\overline{\nu(-y)}f(x-y).
\end{equation*}

The next lemma expresses the autocorrelation of a deterministic measure $\nu\in M_X^{\infty}$ as a limit of orbit averages defined on $M_X^{\infty}$. This result extends \cite[Lem. 2.13]{gouere} to general locally finite atomic measures.

\begin{lemma}\label{g213}
 Let $\gamma$ be a locally finite measure on $\mathbb{R}^d$. Then $\gamma$ is the autocorrelation $\gamma_\nu$ of a (deterministic) measure $\nu\in\mathcal{M}^{\infty}_{X}$ if and only if for all $f\in C_c(\mathbb{R}^d)$,
 
 \begin{equation*}
     \lim_{R\to\infty}\frac{1}{Vol(B_R)}\int_{B_R}H_f(\nu*\delta_{-t})dt=\gamma(f).
 \end{equation*}
\end{lemma}

\begin{proof}
As $\gamma$ is translation bounded, it is sufficient to show that for all $f\in C_c(\mathbb{R}^d)$ (see Lemma 1.2 in \cite{Sch})

\begin{equation}\label{auto-ergod}
    \frac{1}{Vol(B_R)}\left|\int_{B_R}H_f(\nu*\delta_{-x})dt-\sum_{\substack{x\in\mathsf{supp}(\nu)\cap B_R\\ y\in\mathsf{supp}(\nu)}}\nu(x)\overline{\nu(-y)}f(x-y)\right|\longrightarrow 0.
\end{equation}

A simple computation shows that

\begin{equation*}
    \mathcal{H}_f(\nu*\delta_{-t})=\sum_{x,y\in\mathsf{supp}(\nu)}\psi(x-t)\nu(x)\overline{\nu(-y)}f(x-y).
\end{equation*}

On the other hand, for all $R\geq M$ we have the following properties, since $\mathrm{supp}\psi\subset B_M$
\begin{itemize}
    \item If $x\not\in B_{R+M}$, then $\int_{B_R}\psi(x-t)dt=0$, and
    
    \item if $x\in B_{R-M}$, then $\int_{B_R}\psi(x-t)dt=1$.
\end{itemize}

Hence, the left-hand side of (\ref{auto-ergod}) is equal to

\begin{equation*}
\begin{split}
  \frac{1}{Vol(B_R)}\left|\sum_{\substack{x\in\mathsf{supp}(\nu)\cap B_{R+M}\setminus B_{R-M}\\ y\in\mathsf{supp}(\nu)}}f(x-y)\nu(x)\overline{\nu(-y)}\int_{B_{R}}\psi(x-t)dt\right.\\
  -\left.\sum_{\substack{x\in\mathsf{supp}(\nu)\cap B_R\setminus B_{R-M}\\ y\in\mathsf{supp}(\nu)}}\nu(x)\overline{\nu(-y)}f(x-y)\right|.
  \end{split}
\end{equation*}

By triangle inequality, this last expression is bounded by

\begin{equation}\label{auto-ergod-2}
\begin{split}
    \frac{||f||_{\infty}C_{\nu}}{Vol(B_R)}\left(\sum_{x\in\mathsf{supp}(\nu)\cap B_{R+M}\setminus B_{R-M}}\#(\mathsf{supp}(\nu)\cap(\mathsf{supp}(f)+x))\right.\\
    \left.+\sum_{x\in\mathsf{supp}(\nu)\cap B_{R}\setminus B_{R-M}}\#(\mathsf{supp}(\nu)\cap(\mathsf{supp}(f)+x))\right),
\end{split}    
\end{equation}

where $|\nu(x)\overline{\nu(-y)}|\leq C_{\nu}$ is a positive constant only depending on $\nu$.
Moreover since $\nu\in\mathcal{M}^{\infty}_{X}$, $f\in C_c(\mathbb{R}^d)$ and by the discreteness of $X$, it follows that $\#(\mathsf{supp}(\nu)\cap (\mathsf{supp}(f)+x))$ can be bounded by a positive packing constant $C_f$ only depending on $\mathrm{diam}(\mathsf{supp}(f))$ and on the separation constant of $X$. Thus, by writting $D_f=C_f||f||_{\infty}$, we have that (\ref{auto-ergod-2}) can be bounded from above by

\begin{equation}\label{lastex}
    \frac{D_fC_{\nu}}{Vol(B_R)}(\#(\mathsf{supp}(\nu)\cap B_{R+M}\setminus B_{R-M})+\#(\mathsf{supp}(\nu)\cap B_R\setminus B_{R-M})).
\end{equation}

Note that (H.1) implies that for each $M>0$ and for each $p\in\R^d$ there holds
 
\begin{equation}\label{H-4}
    \lim_{R\to\infty}\frac{\# \left(X\cap B_R(p)\setminus B_{R-M}(p)\right)}{Vol(B_R)}=0.
\end{equation}

The proof of \eqref{H-4} is standard: we compare $\#(X\cap B_R(p))$ with the values $\#(X\cap B_{R\pm |p|})$ and compare the latter with $\mathrm{dens}(X) Vol(B_{R\pm|p|})$ via (H.1); then we do the same with $R$ replaced by $R-M$. Due to the precise power law growth $Vol(B_R)=C R^d$ and to (H.1), all the error terms are $o(R^d)$ or $O(R^{d-1})$, which implies \eqref{H-4}.

\medskip

Since $\nu\in\mathcal{M}^{\infty}_{X}$, from property (\ref{H-4}) of $X$ we have that \eqref{lastex} converges to zero when $R$ goes to infinity. Therefore the limit in (\ref{auto-ergod}) holds, as desired. 

\end{proof}

The following result was also proved in \cite{gouere} for point processes, and we note that it extends to random measures, using Lemma 2.13 therein, and it directly extends via Lemma \ref{g213}.

\begin{prop}[{Cf. \cite[Theorem 4.3]{gouere}}]\label{gouere4.3}
Let $\mu$ be a stationary and integrable random measure with values in $\mathcal{M}^{\infty}_X$. Then $\mu$ has almost surely a (random) autocorrelation $\gamma_{\mu}$ such that

\begin{equation*}
    \mathbb{E}(\gamma_{\mu})=I(\widetilde P )\hspace{0.3cm}\mbox{ and }\hspace{0.3cm} \mathbb{E}(\widehat{\gamma_{\mu}})=\widehat{I(\widetilde{P})}.
\end{equation*}
Moreover, if $\mu$ is ergodic, $\gamma_{\mu}=I(\widetilde{P})$ almost surely. In particular, $I(\widetilde{P})$ is positive definite.
\end{prop}

We now proceed with the same strategy as in Section \ref{sec-latcase} to apply Strungaru's criterion, and then conclude as in the beginning of Section \ref{proofs}.

\begin{prop}\label{stationarycase}
Assume the random perturbation $X_{\xi}$ is stationary and weakly-mixing. Then almost surely, for every $\lambda\in\mathbb{R}^d$ the diffraction measure $\widehat{\gamma_{\xi}^{\lambda}}$ has no atoms.
\end{prop}

\begin{proof}
This follows directly from Proposition \ref{gouere4.3} since the stationarity and the weak-mixing condition of $X_{\xi}$ implies that $\mu_{\xi}^{\lambda}$ is a stationary and weakly-mixing random measure for every $\lambda\in\mathbb{R}^d$. The weak-mixing condition of $\mu_{\xi}^{\lambda}$ and Proposition \ref{gouere4.3} implies by expanding the definition like in Lemma \ref{bala_case} that almost surely

\begin{equation*}
    \lim_{R\to\infty}\frac{1}{Vol(B_R)}|\gamma_{\xi}^{\lambda}(B_R)|=0,
\end{equation*}

and therefore by Strungaru's criterion (see \cite[Proposition 4.1]{Strungaru}), we get that $\widehat{\gamma_{\xi}^{\lambda}}$ is almost surely continuous. Finally, the proof that this property holds for every $\lambda\in\mathbb{R}^d$ follows along the same lines as the proof of Proposition \ref{prop2}, since for bounded $f$ we have a control of the $\lambda$-gradients of 

\[
\frac{1}{Vol(B_R)}\int_{B_R}H_f(\mu_{\xi}^{\lambda}*\delta_{-t})dt
\]

as in \eqref{lderbd}. This finishes the proof of Proposition \ref{stationarycase} and of Theorem C.
\end{proof}

%%%%%%%%%%%%%%%%%%%%%%%%%%%%%%%%%%%%%%%%%%%%%%%%%%%%%%%%%%%%%%%%%%%%%%%%%%%%%%%%%%%%%%%%%%%%%%%%%%%%%%%%%%%%%%%%%%%%%%%%%%%%%%%%%%%%%%%%%%%%%%%%%%%%%%%%%%%%%%%%%%%%%%%%%%%%%%%%%%%%%%%%%%%%%%%%%%%%%%%%%%%%%%%%%%%%%%%%%%%%%%%%%%%%%%%%%%%%%%%%%%%%%%%%%%%%%%%%%%%%%%%%%%
\subsubsection{Some examples of stationary random perturbations}\label{exstat}

In this section we provide a class of examples of random perturbations of irreducible cut-and-project sets verifying the hypotheses of Theorem C.\\

Let $L\subset\mathbb{R}^n$ be a lattice and denote by $D_L$ its fundamental domain. Write $\mathbb{R}^n=V_{\mathsf{phys}}\oplus V_{\mathsf{int}}$, where $V_{\mathsf{phys}}\cong\mathbb{R}^d$ and $V_{\mathsf{int}}\cong\mathbb{R}^k$. We denote by $\pi_{\mathsf{phys}}$ and $\pi_{\mathsf{int}}$ the corresponding projections over $\mathbb{R}^d$ and $\mathbb{R}^k$ respectively. Fix $W\subset V_{\mathsf{int}}$ (called in the literature the {\em window}). The {\em cut-and-project} set $X(L,W)\subset V_{\mathsf{phys}}$ is defined by

\begin{equation}\label{cutandproject}
    X=X(L,W):=\pi_{\mathsf{phys}}(L\cap \pi_{\mathsf{int}}^{-1}(W)).
\end{equation}

We say that $X$ is {\em irreducible} if the following conditions hold:
\begin{itemize}
    \item[(CP1) ]$\pi_{\mathsf{int}}(L)$ is dense in $V_{\mathsf{int}}$.
    \item[(CP2) ]$\pi_{\mathsf{phys}}|_L$ is injective.
    \item[(CP3) ]$W$ is a bounded Borel measurable set with non-empty interior and such that its boundary has zero Lebesgue measure in $V_{\mathsf{int}}$. 
\end{itemize}

It is known (see for instance \cite[Theorem 9.4]{BG1}) that under the conditions (CP1)-(CP3) the cut-and-project set (\ref{cutandproject}) is a quasi-periodic set where (H.1) and (H.2) are satisfied automatically. Consider an $L$-stationary random process $(\xi'_q)_{q\in L}$ of identically distributed random vectors defined on $(\Omega,\mathcal{F},\mathbb{P})$ with $\xi'_q\in V_{\mathsf{phys}}$, and let $\eta$ be a random vector independent of $(\xi'_q)_{q\in L}$, which is distributed uniformly over $D_L$. Then $L_{\xi}^{\mathsf{stat}}:=\{q+\xi'_q+\eta:\ q\in L\}$ is $\R^n$-stationary.

\medskip

Set $X_\xi:=\pi_{\mathsf{phys}}(L_\xi^{\mathsf{stat}}\cap (V_{\mathsf{phys}}\times W))$. Note that identifying $V_{\mathsf{phys}}=\R^d$, automatically $X_\xi$ is $\R^d$-stationary, since so is $L_\xi^{\mathsf{stat}}\cap (V_{\mathsf{phys}} \times W)$, as is direct to check. Also note that 
\[
L_\xi^{\mathsf{stat}}\cap (V_{\mathsf{phys}}\times W) = \{q+\xi_q'+\eta:\ q\in L\cap (V_{\mathsf{phys}}\times W-\eta)\}. 
\]
Therefore, in order to have $X_\xi=\{p+\xi_p:\ p\in X\}$ as required for Theorem C, we may set 
\[
\xi_p:=\xi_{\psi_\eta\circ\pi_{\mathsf{phys}}^{-1}(p)}'+\pi_{\mathsf{phys}}(\eta),
\]
in which for each $\eta\in D_L$ the map $\psi_\eta: L\cap (V_{\mathsf{phys}}\times W)\to L\cap (V_{\mathsf{phys}}\times W-\eta)$ is a fixed bijection of bounded distortion. This is possible by assuming that $X$ is up to bounded displacement (BD) equivalent to a lattice, which means that there is a bijection $\phi: X\to \mathcal{L}\subset\mathbb{R}^d$ such that 

\[
\sup_{x\in X}|\phi(x)-x|<\infty.
\]

In fact, if $\phi:X\to\mathcal{L}$ is a BD-equivalence from $X$ to a lattice $\mathcal{L}\subset\mathbb{R}^d$, the conditions (CP1) and (CP3) imply that for every $x\in\mathbb{R}^n$, there is a BD-equivalence $\phi_{x}:\pi_{\mathsf{phys}}(L\cap (V_{\mathsf{phys}}\times W-x))\to\mathcal{L}$ (see \cite[Proposition 2.6]{HKW}). Hence by using (CP2) we can establish the desired bounded distortion bijection $\psi_{\eta}:=\pi_{\mathsf{phys}}^{-1}\circ\phi_{\eta}^{-1}\circ\phi\circ\pi_{\mathsf{phys}}:L\cap (V_{\mathsf{phys}}\times W)\to L\cap (V_{\mathsf{phys}}\times W-\eta)$.

\medskip
If we assume that $\mathbb{E}(|\xi'|^{d+\varepsilon})<\infty$ for some positive number $\varepsilon$, then $\mathbb{E}|\xi|^{d+\varepsilon}<\infty$. Moreover, if $(\xi'_q)_{q\in L}$ is weakly-mixing (with respect to the natural shift in the lattice $L$), then the random perturbation $X_{\xi}$ is weakly-mixing with respect to the $\mathbb{R}^d$-translations. Thus $X_{\xi}$ fulfils the hypotheses of Theorem C and therefore the irreducible cut-and-project set $X$ is almost surely recoverable from $X_{\xi}$.  

%%%%%%%%%%%%%%%%%%%%%%%%%%%%%%%%%%%%%%%%%%%%%%%%%%%%%%%%%%%%%%%%%%%%%%%%%%%%%%%%%%%%%%%%%%%%%%%%%%%%%%%%%%%%%%%%%%%%%%%%%%%%%%%%%%%%%%%%%%%%%%%%%%%%%%%%%%%%%%%%%%%%%%%%%%%%%%%%%%%%%%%%%%%%%%%%%%%%%%%%%%%%%%%%%%%%%%%%%%%%%%%%%%%%%%%%%%%%%%%%%%%%%%%%%%%%%%%%%%%%%%%%

\section{Discussion and open problems}
In this section we discuss the relations of our results with other notions, and we highlight two possible directions of future research.
\subsection{``Quasi-stationarity'' for the recovery of general Meyer sets}\label{stationarygen}
As already described in Section \ref{asyaffine}, recovery Theorems A, B hold also under weaker notions than the independence of the $\xi_p$, and a good weakening is the hypothesis of stationarity and ergodicity. The drawback of this direction is that $X$ needs to have a group structure, at least asymptotically, in order to define stationarity. 

\medskip 
It seems natural to look for a definition of stationarity of the $(\xi_p)_{p\in X}$ that together with an associated notion of ergodicity, allows to obtain recovery Theorems A, B for more general quasi-periodic $X$, beyond asymptotically affine deformations of lattices. More precisely, if $X$ is quasi-periodic then we may look for notions of invariance of the $(\xi_p)_{p\in X}$ in terms of the geometry of the almost-periods of $X$, which we could refer to as "quasi-stationarity" (however it seems to us that this term is already in use with a different meaning). Perhaps this can allow a more explicit recovery result for "quasi-stationary" $(\xi_p)_{p\in X}$, as opposed to Theorem C, in which stationarity under the $\R^d$-action is imposed on the perturbed $X_\xi$, and the relation between the geometry of $X$ and that of the $(\xi_p)_{p\in X}$ is not explicitly quantified.

\medskip
A natural setup where to start is that $X$ satisfies the hypotheses (H.1) and (H.2), and that $X$ is a {\em Meyer set}, which means that $X$ is relatively dense and $X-X:=\{p-q:\ p,q\in X\}$ is uniformly discrete. %Observe that a Meyer set satisfies automatically condition (H.2) of Section \ref{ssec:hypo} (see for instance Theorem 4.1 in \cite{Lag}).\\

If $X-X$ is uniformly discrete then a first simplification occurs, since the limit defining $\gamma_\xi^\lambda(k)$ is trivially zero for $k\notin X-X$. The difficulty is that, in general, for $k\in X-X$ one can not ensure the existence and good control of the limits
\begin{equation}\label{autocorr}
\gamma_\xi^\lambda(k)=\lim_{R\to \infty} \frac{1}{Vol(B_R)} \sum_{\substack{p\in X\cap B_R\\ p-k\in X}}\mu_\xi^\lambda(p)\overline{\mu_\xi^\lambda(p-k)}.
\end{equation}
In this case, the main open question is to find explicit further conditions that we need to impose on $X$ and on the $\xi_p$ in order to achieve the result of Proposition \ref{prop2}. 

\subsection{Recovering the structure factor and almost-sure recovery of the diffraction measure}\label{recoverdiffraction}

In the setting of Theorem A, if the characteristic function $\varphi(\lambda):=\mathbb E[e^{-2\pi i\langle \xi, \lambda\rangle}]$ vanishes at some points of the support of $\widehat \delta_X$, this means that this part of the spectrum of $X$ is not recoverable, a phenomenon called \emph{cloaking} of $X$ by the random perturbations $\xi_p$.

\medskip
In \cite{kkt} the authors present a result analogous to our theorems, but in \emph{expectation form}, formulated in terms of \emph{structure factors}. In our notation, the structure factor of a random discrete set $X_\xi=\{p+\xi_p:\ p\in X\}$ is the function $S_{X_\xi}:\mathbb R^d\to[0,+\infty]$ formally given by
\begin{equation}\label{strfact}
S_{X_\xi}(\lambda):=\lim_{R\to \infty}\mathbb E\left[\frac{1}{\#(X\cap B_R)}\left|\sum_{p\in X\cap B_R}e^{-2\pi i \langle p+\xi_p, \lambda\rangle}\right|^2\right]=\frac{\mathbb E[\widehat{\gamma}_{X_\xi}(\lambda)]}{\mathrm{dens}(X)},
\end{equation}
where $\gamma_{X_\xi}(\lambda)$ is the autocorrelation of $X_{\xi}$. Note that the above limit may not exist, or may be defined only in a weak (distributional) sense, if $X_\xi$ is not regular enough.
\medskip

The above quantity is most commonly studied for the case that $X$ is not deterministic, but has a well-defined density, and $\mathrm{dens}(X)=\mathrm{dens}(X_\xi)=\widehat\gamma_{X_\xi}(0)$ almost surely. Often $X$ is assumed to be a stationary spatial point process, the simplest example being a "stationarized lattice", i.e. $X$ is lattice translated by a random vector $\tau$ uniformly distributed over its fundamental domain. In this case, if the $\xi_p$ are assumed to be i.i.d. and independent of (the random process) $X$, then as done in \cite[App. A]{kkt} one can rewrite the square modulus from \eqref{strfact} as a double sum, and then split the expectation as a product of expectations expressed in terms of $X, \xi_p, p\in X$. Thus we get directly (where the $-1$ terms come from the diagonal term in the double sums obtained expanding the square in \eqref{strfact})
\begin{equation}\label{strfact2}
S_{X_\xi}(\lambda) - 1 = \left|\mathbb E\left[e^{-2\pi i\langle \xi, \lambda\rangle}\right]\right|^2 \left(S_X(\lambda) - 1\right).
\end{equation}
Our Theorems A, B, C give conditions for the almost-sure recovery of $\widehat{\delta_X}$, the Fourier transform of the empirical measure of the set $X$. The same strategy does not directly give the almost-sure recovery of the diffraction measure of $X$, which would be a strengthening of \eqref{strfact2}, and it would be interesting to continue our study in that direction. More precisely, this requires to establish conditions on $X,(\xi_p)_{p\in X}$ such that a limit like the below is ensured almost surely for all $\lambda\in\R^d$:
\begin{equation*}
    \lim_{R\to\infty}\frac{1}{|B_R|}\sum_{\substack{p,q\in X\cap B_R}}e^{-2\pi i\langle p-q,\lambda\rangle}\left(e^{-2\pi i\langle\xi_p-\xi_q,\lambda\rangle}-\left|\mathbb{E}[e^{-2\pi i\langle\xi,\lambda\rangle}]\right|^2\right)=0.
\end{equation*}
If we prove the above, then this would substitute Proposition \ref{reduction Main Thm}, and would imply by the same reasoning as in Section \ref{proofs} that almost surely for all $\lambda$

\begin{equation}\label{recovdiff}
\widehat{\gamma_{X_\xi}}(\lambda) - \mathsf{dens}(X)= \left|\mathbb E\left[e^{-2\pi i\langle \xi, \lambda\rangle}\right]\right|^2 \left(\widehat{\gamma_X}(\lambda)-\mathsf{dens}(X)\right),
\end{equation}

which is the almost-sure version of \eqref{strfact2} and which already appears in \cite{Hof} under different hypotheses over the set $X$. In fact, in \cite{Hof} it is assumed that $X$ is an uniformly discrete set with asymptotic density and satisfying an ergodicity condition. Hof's ergodicity hypotheses are not directly comparable to our (H1), (H2). We leave to future work the elucidation of implications between hypotheses of ergodicity of $\delta_X$ like those in \cite{Hof} and hypotheses on $\widehat{\gamma}_X$ such as (H2).
\medskip

%%%%%%%%%%%%%%%%%%%%%%%%%%%%%%%%%%%%%%%%%%%%%%%%%%%%%%%%%%%%%%%%%%%%%%%%%%%%%%%%%%%%%%%%%%%%%%%%%%%%%%%%%%%%%%%%%%%%%%%%%%%%%%%%%%%%%%%%

%%%%%%%%%%%%%%%%%%%%%%%%%%%%%%%%%%%%%%%%%%%%%%%%%%%%%%%%%%%%%%%%%%%%%%%%%%%%%%%%%%%%%%%%%%%%%%%%%%%%%%%%%%%%%%%%%%%%%%%%%%%%%%%%%%%%%%%%%%%%%%%%%%%%%%%%%%%%%%%%%%%%%%%%%%%%%%%%%%%%%%%%%%%%%%%%%%%%%%

\appendix

\section{Auxiliary results}
\subsection{Strong law of large numbers}

\begin{prop}\label{prop:slln}
Let $(X_n)_{n\in\mathbb{N}}$ be a sequence of real-valued random variables such that almost surely $X_n\ge C_n$ for suitable constants $C_n\in\R$. Let $S_n:=\sum_{k=1}^n (X_k-C_k)$ and assume that the following hold.
\begin{enumerate}
    \item $\sup_{n\in\N}\frac1n\mathbb{E}(S_n)<\infty$.
    \item For any $\alpha>1$ if $k_n:=[\alpha^n]$, where the brackets $[\cdot]$ denote the integer part of a real number, there holds
    \[
    \sum_{n=1}^\infty\frac{\mathsf{Var}(S_{k_n})}{k_n^2}<\infty.
    \]
\end{enumerate}
Then almost surely
\begin{equation}\label{chgos}
\frac1n \sum_{k=1}^n\left(X_k - \mathbb E[X_k]\right) \to 0.
\end{equation}
\end{prop}

\begin{proof}
Without loss of generality assume that $C_n=0$ for every $n\in\mathbb{N}$. For any positive number $\varepsilon$, Chebyshev's inequality yields

\begin{equation*}\label{var-cov}
\begin{split}
    \sum_{n\geq 1}\mathbb{P}\left(\left|\dfrac{S_{k_n}-\mathbb{E}(S_{k_n})}{k_n}\right|>\varepsilon\right)&\leq \frac{1}{\epsilon^2}\sum_{n\geq 1}\dfrac{\mathsf{Var}(S_{k_n})}{k_n^2}<\infty.
    \end{split}
\end{equation*}
The rest of the proof follows the same steps of \cite[Thm. 1]{Ch-Gos} or \cite{Et}.
\end{proof}

\begin{prop}\label{Etem-Cov}
Let $(X_n)_{n\geq 1}$ be a sequence of identically distributed positive random variables, and denote $S_n=\sum_{i=1}^n X_i$. Moreover assume that $\mathbb{E}[X_1]<\infty$ and for $Y_i:=X_i\mathbbm{1}_{\{X_i\leq i\}}$ suppose that
\begin{equation*}
    \sum_{n\ge 1}\frac1{n^2}\sum_{\substack{i,j=1\\i\neq j}}^n|\mathsf{Cov}(Y_i,Y_j)|<\infty.
\end{equation*}
Then $S_n/n$ converges almost surely to $\mathbb{E}[X_1]$.
\end{prop}

\begin{proof}
Define $S_n^*:=\sum_{i=1}^nY_i$. For $k_n$ as in Proposition \ref{prop:slln} we have that

\begin{equation*}\label{var-cov2}
\begin{split}
    \sum_{n\geq 1}\mathbb{P}\left(\left|\dfrac{S^*_{k_n}-\mathbb{E}(S^*_{k_n})}{k_n}\right|>\varepsilon\right)&\leq \frac{1}{\epsilon^2}\sum_{n\geq 1}\dfrac{\mathsf{Var}(S_{k_n}^*)}{k_n^2}\\
    &=\frac{1}{\epsilon^2}\sum_{n\geq 1}\dfrac{1}{k_n^2}\left(\sum_{i=1}^{k_n}\mathsf{Var}(Y_i)+\sum_{\substack{i,j=1\\ i\neq j}}^{k_n} \mathsf{Cov}(Y_i,Y_j)\right).
    \end{split}
\end{equation*}

By the hypotheses on $(X_n)_{n\geq 1}$, we have that 
\begin{equation}\label{etem-proof}
    \sum_{n\geq 1}\frac{1}{k_n^2} \sum_{\substack{i,j=1\\ i\neq j}}^{k_n}|\mathsf{Cov}(Y_i,Y_j)|<\infty,
\end{equation}
which implies the convergence of the second sum in (\ref{etem-proof}). The rest of the proof follows the very same lines than Theorem 1 in \cite{Et}.
\end{proof}

\begin{remark}
The hypotheses of Proposition \ref{Etem-Cov} are satisfied if $(X_n)_{n\geq 1}$ are positive, identically distributed with correlation $|\mathsf{Cov}(Y_i,Y_j)|\le \beta^{|i-j|}$ for some $0<\beta<1$, and if $\mathbb{E}[X_1]<\infty$.
\end{remark}
\subsection{Hellinger density of diffraction measures bounds }

\medskip

Let $\gamma,\gamma'$ be locally finite measures. We define the localized Hellinger density $\rho(\gamma,\gamma')$, which is a locally finite positive measure, by requiring, for any locally finite measure $\sigma$ over $\R^d$ such that $\gamma\ll\sigma $ and $\gamma'\ll \sigma$,
\begin{equation}
    \forall f\in \mathcal{S}(\R^d),\quad \langle \rho(\gamma,\gamma'),f\rangle:=\left\langle \left(\frac{d\gamma}{d\sigma}\right)^{1/2}\left(\frac{d\gamma'}{d\sigma}\right)^{1/2}d\sigma, f\right\rangle.
\end{equation}
Note that if $\gamma,\gamma'$ are mutually singular, then directly from the definition it follows that $\rho(\gamma,\gamma')=0$.

\medskip

We claim that the following holds
\begin{thm}\label{yakirlemma4}
 If $\mu,\nu$ are translation-bounded, Fourier-transformable measures over $\R^d$ and $\widehat{\gamma_\mu},\widehat{\gamma_\nu}$ are their diffraction measures, then 
 \begin{equation}
     \forall f \in C_c(\R^d), f\ge 0\quad \limsup_{R\to\infty}\left|\frac{1}{Vol(B_R)}\iint_{B_R\times B_R}\widehat f(p-q)\ d\mu(p)\ d\overline{\nu(q)}\right|\le \langle \rho(\widehat{\gamma_\mu},\widehat{\gamma_\nu}), f\rangle.
 \end{equation}
For $\mu,\nu\ll \sum_{p\in X}\delta_p$ (i.e. $\mu,\nu$ atomic with atoms $\subset X$):
 \begin{equation}\label{helldens2}
     \forall f \in \mathcal{S}(\R^d)\quad \limsup_{R\to\infty}\left|\frac{1}{Vol(B_R)}\sum_{p,q\in X\cap B_R}\mu(p)\overline{\nu(p)} \widehat f(p-q)\right|\le \langle \rho(\widehat{\gamma_\mu},\widehat{\gamma_\nu}), |f|\rangle.
 \end{equation}
\end{thm}
As we prove before, Proposition \ref{reduction Main Thm} follows directly from theorem \ref{yakirlemma4}, applied to $\mu_X^\lambda, \mu_\xi^\lambda$, once we know that the $\widehat{\gamma_X^\lambda}$ is purely atomic and that $\widehat{\gamma_\xi^\lambda}$ is absolutely continuous with respect to Lebesgue measure.

\medskip
We use the below lemma which follows like Theorem B from Yakir:\\

\begin{lemma}\label{yakirlthmb}
 Let $\gamma_R,\gamma_R'$ be two families of positive measures over $\R^d$ such that $\gamma_R\to\gamma$ and $\gamma_R'\to \gamma'$ vaguely as $R\to\infty$. Then
 \begin{equation*}
     \forall f\in \mathcal{S}(\R^d),  f\ge 0, \quad \limsup_{R\to\infty}\langle \rho(\gamma_R,\gamma_R'), f\rangle\le \langle \rho(\gamma,\gamma'), f\rangle. 
 \end{equation*}
\end{lemma}
\medskip

\begin{proof}
 Consider the locally finite measure $\sigma_R:=\gamma_R+\gamma_R'$. It is direct that $\gamma_R$ and $\gamma_R'$ are absolutely continuous with respect to $\sigma_R$. Thus, from the Cauchy-Schwartz inequality we get
 
 \begin{equation}\label{limsup_hell}
     \begin{split}
         \langle\rho(\gamma_R,\gamma_R'),f\rangle &=\displaystyle\int f(x)\left(\dfrac{d\gamma_R}{d\sigma_R}\right)^{1/2}\left(\dfrac{d\gamma_R'}{d\sigma_R}\right)^{1/2}d\sigma_R\\
         &\leq \left( \displaystyle\int f(x)\left(\dfrac{d\gamma_R}{d\sigma_R}\right)d\sigma_R\right)^{1/2}\left(\int f(x)\left(\dfrac{d\gamma_R'}{d\sigma_R}\right)d\sigma_R\right)^{1/2}\\
         &=\left( \displaystyle\int f(x)d\gamma_R\right)^{1/2}\left( \displaystyle\int f(x)d\gamma_R'\right)^{1/2}
         \end{split}
 \end{equation}

On the other hand, define:

\begin{equation*}
\begin{split}
    A :=\left\{x\in\mathbb{R}^d:\ \frac{d\gamma}{d\sigma}(x)=0\right\},\hspace{1cm}  B:=\left\{x\in\mathbb{R}^d:\ \frac{d\gamma'}{d\sigma}(x)=0\right\}\setminus A,
    \end{split}
\end{equation*}
and for $\varepsilon>0$ we consider

\begin{equation*}
    V_j(\varepsilon):=\left\{x\in\mathbb{R}^d:\ (1+\varepsilon)^{j-1}\dfrac{d\gamma}{d\sigma}(x)\leq\dfrac{d\gamma'}{d\sigma}(x)\leq(1+\varepsilon)^{j}\dfrac{d\gamma}{d\sigma}(x)\right\}\setminus (A\cup B).
\end{equation*}

and $f_j:=f\mathbbm{1}_{V_j}$. Note that the family $$\{A, B,\{V_j(\varepsilon)\}_{j\in\mathbb{Z}}\}$$ determines a partition of $\mathbb{R}^d$.
From (\ref{limsup_hell}) and the definition of vague convergence we have that for each $f_j$ 
 
 \begin{eqnarray*}
\lefteqn{\limsup_{R\to\infty}\langle\rho(\gamma_R,\gamma_R'),f_j\rangle\leq\left( \displaystyle\int f_j(x)d\gamma\right)^{1/2}\left( \displaystyle\int f_j(x)d\gamma'\right)^{1/2}}\\
     &=&\left( \displaystyle\int_{V_j} f(x)\dfrac{d\gamma}{d\sigma}d\sigma\right)^{1/2}\left( \displaystyle\int_{V_j} f(x)\dfrac{d\gamma'}{d\sigma}d\sigma\right)^{1/2}\\
     &\leq&\left( \displaystyle\int_{V_j} (1+\varepsilon)^{-(j-1)/2}f(x)\left(\dfrac{d\gamma}{d\sigma}\right)^{1/2}\left(\dfrac{d\gamma'}{d\sigma}\right)^{1/2}d\sigma\right)^{1/2}\left( \displaystyle\int_{V_j} (1+\varepsilon)^{\frac{j}{2}}f(x)\left(\dfrac{d\gamma'}{d\sigma}\right)^{1/2}\left(\dfrac{d\gamma}{d\sigma}\right)^{1/2}d\sigma\right)^{1/2}\\
     &=&(1+\varepsilon)^{1/4}\displaystyle\int f_j(x)\left(\dfrac{d\gamma}{d\sigma}\right)^{1/2}\left(\dfrac{d\gamma'}{d\sigma}\right)^{1/2}d\sigma\\
     &=&(1+\varepsilon)^{1/4}\langle\rho(\gamma,\gamma'),f_j\rangle.
 \end{eqnarray*}
Furthermore, for $f_A:=f\mathbbm{1}_{A}$ and $f_B:=f\mathbbm{1}_B$ we obtain

\begin{equation*}
    \langle\rho(\gamma_R,\gamma'_R),f_A\rangle=\langle\rho(\gamma_R,\gamma'_R),f_B\rangle=0.
\end{equation*}

Therefore, by summing over $j$ and using the fact that the sum of the limsup's is larger or equal than the limsup of the sum, and then taking $\epsilon\downarrow 0$, the thesis follows.
\end{proof}

\begin{proof}[Proof of Theorem \ref{yakirlemma4} from Lemma \ref{yakirlthmb}:]
 As in \cite{Yakir}, it suffices to consider a sequence of measures $\gamma_{\mu}^{R}$ and $\gamma_{\nu}^R$ absolutely continuous with respect the Lebesgue measure whose Fourier transform converge vaguely to the diffraction measures $\widehat{\gamma_{\mu}}$, $\widehat{\gamma_{\nu}}$, respectively.\\
 
 The key point here is the fact that the Fourier transform of a finite pure point measure is proportional to the Lebesgue measure. This, together the identity (\ref{F-lim=lim-F}) yield 
 
 \begin{equation*}
     \langle\rho(\widehat{\gamma_{\mu}^{R}},\widehat{
     \gamma_{\mu}^{R}}),f\rangle\geq  \left|\dfrac{1}{Vol(B_R)}\sum_{p,q\in X\cap B_R}\mu(p)\overline{\nu(q)}\widehat{f}(p-q)\right|.
 \end{equation*}
 
 In fact, let $\mu, \nu$ be two atomic measures with atoms in $X$. Write 
 
 \begin{equation*}
 \mu_R =\sum_{p\in X\cap B_R}\mu (p)\delta_p,\qquad \nu_R=\sum_{p\in X\cap B_R}\nu (p)\delta_p.
 \end{equation*}
 
 Thus, from (\ref{F-lim=lim-F}), we have that the diffraction measures $\widehat{\gamma_{\mu}}$ and $\widehat{\gamma_{\nu}}$ are vague limits as follows
 
 \begin{equation*}
     \widehat{\gamma_{\mu}}=\lim_{R\to \infty}\dfrac{1}{Vol(B_R)}\left|\sum_{p\in X\cap B_R}\mu(p)e^{-2\pi i\langle p,x\rangle}\right|^2 dx,\quad
     \widehat{\gamma_{\nu}}=\lim_{R\to \infty}\dfrac{1}{Vol(B_R)}\left|\sum_{p\in X\cap B_R}\nu(p)e^{-2\pi i\langle p,x\rangle}\right|^2 dx.
 \end{equation*}
 
 Then, for every $f\in\mathcal{S}(\mathbb{R}^d)$, $f\geq 0$ 
 
 \begin{equation*}
     \begin{split}
         \langle\rho(\widehat{\gamma_{\mu}^R}, \widehat{\gamma_{\nu}^R}),f\rangle &=\dfrac{1}{Vol(B_R)}\displaystyle\int f(x)\left|\sum_{p\in X\cap B_R} \mu(p)e^{-2\pi i\langle p,x\rangle}\right|\left|\sum_{q\in X\cap B_R} \nu(q)e^{-2\pi i\langle q,x\rangle}\right|dx\\
         &\geq \left|\dfrac{1}{Vol(B_R)}\sum_{p,q\in X\cap B_R}\mu(p)\overline{\nu(q)}\displaystyle\int f(x)e^{-2\pi i\langle x,p-q\rangle} dx\right|\\
         &= \left|\dfrac{1}{Vol(B_R)}\sum_{p,q\in X\cap B_R}\mu(p)\overline{\nu(q)}\widehat{f}(p-q)\right|.
    \end{split}
 \end{equation*}

Hence by taking $\limsup$ over $R$ and by Lemma \ref{yakirlthmb}
we get (\ref{helldens2}) as claimed.
\end{proof}

\end{document}